\documentclass[smallcondensed]{svjour3}
\usepackage[utf8]{inputenc}

\usepackage{amsmath,amssymb,amsfonts,mathtools,float,mathdots}
\usepackage{algorithm}
\usepackage{appendix}
\usepackage{algcompatible}

\newcommand{\mc}[1]{\mathcal{#1}}
\newcommand{\mb}[1]{\mathbb{#1}}
\newcommand{\bmat}[1]{\begin{bmatrix} #1 \end{bmatrix}}
\renewcommand{\vec}[1]{\mathbf{#1}}
\newcommand{\mat}[1]{\mathbf{#1}}
\newcommand{\ten}[1]{\boldsymbol{\mathcal{#1}}}
\newcommand{\matc}[1]{\boldsymbol{\mathcal{#1}}}
\newcommand{\tenh}[1]{\widehat{\boldsymbol{\mc{#1}}}}
\newcommand{\rank}{\mathsf{rank}}

\newcommand{\mathat}[1]{\mathbf{\widehat{#1}}}

\begin{document}

\title{Efficient randomized tensor-based algorithms for function approximation and low-rank kernel interactions\thanks{A.K.S. was supported, in part, by the National Science Foundation (NSF) through the grant DMS-1821149 and DMS-1745654. M.E.K. was supported, in part, by the National Science Foundation under Grants NSF DMS 1821148 and Tufts T-Tripods Institute NSF HDR grant CCF-1934553.}}
\author{Arvind K. Saibaba \and Rachel Minster \and   Misha E.\ Kilmer}
\institute{
Arvind K.\ Saibaba \at Department of Mathematics \\ North Carolina State University, USA\\\email{asaibab@ncsu.edu} \\ORCID: 0000-0002-8698-6100
\and
Rachel Minster \at Department of Computer Science \\ Wake Forest University, USA \\ \email{minsterr@wfu.edu} 
\and
Misha E.\ Kilmer \at Department of Mathematics \\  Tufts University, USA \\ \email{misha.kilmer@tufts.edu}
} 
\titlerunning{Randomized tensor-based algorithms function and kernel approximations}
\date{Received: date / Accepted: date}

\maketitle

\begin{abstract}
In this paper, we introduce a method for multivariate function approximation using function evaluations, Chebyshev polynomials, and tensor-based compression techniques via the Tucker format.  We develop novel randomized techniques to accomplish the tensor compression, provide a detailed analysis of the computational costs, provide insight into the error of the resulting approximations, and discuss the benefits of the proposed approaches. We also apply the tensor-based function approximation to develop low-rank matrix approximations to kernel matrices that describe pairwise interactions between two sets of points; the resulting low-rank approximations are efficient to compute and store (the complexity is linear in the number of points). We have detailed numerical experiments on example problems involving multivariate function approximation, low-rank matrix approximations of kernel matrices involving well-separated clusters of sources and target points, and a global low-rank approximation of kernel matrices with an application to Gaussian processes. 
    

\end{abstract}

\section{Introduction}
Multivariate function approximation to construct a surrogate model for an underlying physical model that is expensive to evaluate is an important topic: it has applications in optimization, reduced order modeling, uncertainty quantification, sensitivity analysis, and many other areas.  The first goal of this paper is to provide a multivariate approximation of a function $f: \mc{D} \times \mb{R}$ of the form 
\begin{equation}\label{eqn:multi_approx} f(x_1,\dots,x_N) \approx \sum_{i_1=1}^{r_1} \dots \sum_{i_N=1}^{r_N} c_{i_1,\dots,i_N} \phi_{i_1}^{(1)}(x_1) \cdots  \phi_{i_N}^{(N)}(x_N),\end{equation}
where $\mc{D} \subset \mb{R}^N$ is a hypercube, $\{\phi_{j}^{(k)}\}$ are basis functions for $j=1,\dots,r_k$ and $k=1,\dots,N$ and $c_{i_1,\dots,i_N}$ are the coefficients. Such a decomposition is described as a summation of separable functions and the existence of such an approximation depends on the properties of $f$ (such as regularity). If the functions are analytic, then the error in the best approximation decays geometrically with the degree of the polynomial~\cite{trefethen2017multivariate,gass2018chebyshev}; the decay is algebraic if the function is in a Sobolev space~\cite{griebel2019analysis}.

This paper proposes a technique for multivariate function approximation, using a combination of Chebyshev polynomial interpolation using function evaluations at Chebyshev interpolation nodes  and compression of the tensor (i.e. multiway array) that contains the function evaluations.  We use the Tucker format for the tensor to achieve compression. Since it is a natural fit for the application at hand, we can leverage prior work on randomized low-rank tensor approximations~\cite{minster2020randomized} and we can exploit the structure in the Tucker format to our computational advantage resulting in new algorithms. Similar functional tensor approximations have been used for trivariate functions (i.e., for $N=3$)~\cite{dolgov2021functional,hashemi2017chebfun} but our approach is not limited to three dimensions. The authors in~\cite{rai2019randomized} also use a functional Tucker approximation but use sparsity promoting techniques to reduce the storage cost involving the core tensor. Several authors have constructed multivariate tensor approximations using the Tensor Train format, rather than the Tucker format~\cite{gorodetsky2019continuous,bigoni2016spectral}; however, a detailed comparison between the Tucker and Tensor Train formats is beyond the scope of the paper. In addition to tensor-based methods, there is a technique for multivariate function approximation that uses sparse grid interpolation~\cite{bungartz_griebel_2004} but it makes slightly different assumptions on the regularity of the functions to be approximated. While there are alternatives to the Tucker format for function approximation, to the best of our knowledge, none of the tensor formats have been explored in the context of kernel approximations, which is the second goal of the paper. For this application, using the Tucker format for the tensor has several advantages, which we exploit throughout this paper.

Kernel methods are used to model pairwise interactions among a set of points. They are popular in many applications in numerical analysis and scientific computing such as 
solvers for integral equations, radial basis function interpolation, Gaussian processes, machine learning, geostatistics and spatiotemporal statistics, and Bayesian inverse problems. A major challenge in working with kernel methods is that the resulting kernel matrices are dense, and are difficult to store and compute with when the number of interaction points is large. Many fast kernel summation methods have been developed over the years including the Barnes-Hut algorithm \cite{barnes1986hierarchical}, Fast Multipole Method \cite{greengard1987fast}, Fast Gauss Transform \cite{greengard1991fast}, etc.

 Operations with kernel matrices, which are used to model pairwise interactions, can be accelerated by using low-rank matrix approximations.  
 Over the years, there have been many efforts to compute efficient low-rank approximations in the context of kernel methods~\cite{cambier2019fast,xu2018low,chen2018fast,ying2004kernel,ye2020analytical,bebendorf2009recompression}.
 There are two techniques to generating such approximations.  In the first approach, the entire matrix is approximated by a global low-rank matrix (e.g., using Nystr\"om approach). In the second approach,  the kernel matrix is represented as a hierarchy of subblocks, and certain of the off-diagonal subblocks are then approximated in a low-rank format (e.g., \cite{hackbusch1999sparse,hackbusch2002data,grasedyck2003construction,borm2003hierarchical,hackbusch2015hierarchical}).  The subblocks representing the interactions between these sets of points can be compressed using the singular value decomposition (SVD). However, the cost of computing the SVD can be prohibitively expensive when the number of source and target points is large and/or when the low-rank approximations need to be computed repeatedly over several sets of source and target points.  Our goal is to avoid these expensive calculations when computing our approximations.

Our approach builds on the black-box fast multipole method (BBFMM) low-rank approximation approach in Fong and Darve~\cite{fong2009black}, which uses Chebyshev interpolation to approximate the kernel interactions. In this method, instead of computing the kernel interactions between every pair of points given, the interactions are only computed between points on Chebyshev grids. This significantly reduces the number of kernel evaluations needed. This work is also similar to the semi-analytic methods in~\cite{cambier2019fast,xu2018low} which also use some form of Chebyshev interpolation to reduce the number of kernel interactions. Our work differs in the subsequent low-rank approximation techniques, as we use tensor-based methods instead of the matrix techniques presented in these papers.

\paragraph{Contributions and Contents.} This paper provides new algorithms for multivariate interpolation with an application to low-rank kernel matrix approximations. We summarize the main contributions in our paper:

\begin{enumerate}
    \item We propose a method for multivariate function interpolation that combines Chebyshev polynomials and tensor compression techniques via the Tucker format. We analyze the error in the approximation which is comprised of the error due to interpolation and the error due to tensor compression. The resulting function approximations can be used as surrogates in a variety of applications including, but not limited to, uncertainty quantification and reduced order modeling. 
    \item We develop three novel randomized compression algorithms for computing low multirank approximations to tensors in the Tucker format. The first algorithm  employs randomized interpolatory decomposition for each mode-unfolding;  we derive error bounds, in expectation, for this compression technique. The second algorithm uses the first approach but for a subsampled tensor, and requires fewer function evaluations. The third algorithm is similar to the first but uses random matrices formed using Kronecker products, and requires less storage and fewer random numbers to be generated. We also provide a detailed analysis of the computational costs. {\it These algorithms are of independent interest to tensor compression beyond computing functional approximations. }
    \item  We use the multivariate function approximations to develop low-rank approximations to kernel matrices which represent pairwise interactions. The resulting approximations are accurate, as demonstrated by detailed numerical experiments, and the computational costs are linear in the number of source and target points. To our knowledge, this is the first use of tensor-based function approximations to low-rank kernel matrix approximations.
\end{enumerate}
We conclude this section with a brief organizational overview of the paper. In Section~\ref{sec:back}, we review the material on Chebyshev interpolation, tensor compression, and randomized matrix algorithms. In Section~\ref{sec:lowrankalgs}, we propose a tensor-based method for approximating multivariate functions, which combines multivariate Chebyshev interpolation with randomized compression algorithms. We develop three randomized algorithms for compressing low-rank tensors and provide insight into the error in the interpolation. In Section~\ref{sec:kernel}, we apply the tensor-based algorithms developed in Section~\ref{sec:lowrankalgs} to develop low-rank approximations to kernel matrices. We provide a detailed analysis of the computational cost of the resulting algorithms. In Section~\ref{sec:results}, we perform three sets of numerical experiments: the first highlights the performance of the tensor-based function approximation, the second demonstrates the accuracy of the low-rank approximations to kernel matrices with several different kernels, and the third is an application to Gaussian processes involving low-rank kernel matrices.

\section{Background}\label{sec:back}
In this section, we provide some background on Chebyshev interpolation as well as tensor notation and decompositions. We also discuss the randomized matrix algorithms that will be important components in our algorithms.
\subsection{Chebyshev interpolation}We begin by reviewing how to interpolate a function $g$: $[a,b] \to \mb{R}$ using Chebyshev interpolation. The first step is to construct the Chebyshev points of the first kind in the interval $[-1,1]$; these are given by $\xi_k =  \cos\left( \frac{2k-1}{2n}\pi\right) $ for $k=1,\dots,n$. Using the mapping $I_{[a,b]}$: $[-1,1] \to [a,b]$ defined as 
\begin{equation}\label{eq:Imap}
I_{[a,b]}(x) = (x+1)\frac{(b-a)}{2}+a,
\end{equation}
we obtain the Chebyshev points in the interval $[a,b]$ as $\eta_k = I_{[a,b]}(\xi_k)$ for $k=1,\dots,n$. This mapping is invertible, with the inverse $I_{[a,b]}^{-1}(x) = 2\frac{(x-a)}{(b-a)}-1$. Recall that for a function $g(x)$ with $x \in [a,b]$, the Chebyshev interpolation polynomial of degree $n-1$ is 
$$\pi_{n-1}(x) = \sum_{k=0}^{n-1} c_k T_k(I^{-1}_{[a,b]}(x)),$$
where $T_k(x) = \cos(k\cos^{-1}(x))$ is the Chebyshev polynomial of degree $k$, and $c_k$ are the coefficients of the Chebyshev polynomials. The coefficients $c_k$ can be obtained using the interpolating conditions and the discrete orthogonality of the Chebyshev polynomials. 
 If we define an interpolating polynomial $S_n^{[a,b]}(x,y)$ as
\begin{equation}\label{eqn:Sn}
    S_n^{[a,b]}(x,y) = \frac{1}{n} + \frac{2}{n} \sum_{k=1}^{n-1} T_k\left(I^{-1}_{[a,b]}(x)\right) \, T_k\left(I^{-1}_{[a,b]}(y)\right),
\end{equation}
then the polynomial $\pi_{n-1}(x)$ can also be expressed as 
\begin{equation}\label{eq:1dinterp}
    \pi_{n-1}(x) = \sum_{k=1}^n g(\eta_k) S_n^{[a,b]}(\eta_k,x).
\end{equation}
The idea of interpolation using Chebyshev polynomials can be extented to multivariate functions as we explore in Section~\ref{ssec:multi_cheb}. The accuracy of Chebyshev interpolation has been studied extensively; see for example~\cite{mason2002chebyshev,trefethen2013approximation}.  
\subsection{Tensor Preliminaries}
We now introduce notation and background information for handling tensors. See \cite{kolda2009tensor} for more details. A $N$-mode tensor is denoted as $\ten{X} \in \mb{R}^{I_1 \times I_2 \times \dots \times I_N}$ with entries $x_{i_1 \dots i_N}$ for $1 \leq i_j \leq I_j$, $1\leq j \leq N$. 

\paragraph{Matricization} We can reorder the elements of a tensor to ``unfold'' it into a matrix, and this process is known as matricization or mode-unfolding. A $d$-mode tensor (also called a $d$th order tensor) can be unfolded into a matrix in $d$ different ways. For mode-$j$, the mode-$j$ fibers\footnote{A mode-$j$ fiber of a tensor is obtained by holding all indices except the $j$th fixed. For example, the mode-1 fibers of 3rd order tensor $\ten{X}$ would be denoted 
$\ten{X}_{:,i_2,i_3}$.} of a tensor $\ten{X}$ are arranged to become the columns of a matrix. This matrix is called the mode-$j$ unfolding, and is denoted by $\mat{X}_{(j)}\in \mb{R}^{I_j \times \prod_{k \neq j} I_k}$. We also denote the unfolding of modes $1$ through $j$ as $\mat{X}_{(1:j)} \in \mb{R}^{ (I_1\dots I_j) \times (I_{j+1}\dots I_N)}$.

\paragraph{Tensor product} One fundamental operation with tensors is the mode-wise product, a way of multiplying a tensor by a matrix. For a matrix $\mat{A} \in \mb{R}^{m \times I_j}$, the mode-$j$ product of $\ten{X}$ with $\mat{A}$ is denoted $\ten{Y} = \ten{X} \times_j \mat{A} \in \mb{R}^{I_1 \times \dots \times I_{j-1} \times m \times I_{j+1} \times \dots \times I_N}$. The entries of this product are 
$$\ten{Y}_{i_1,\dots,i_{j-1},k,i_{j+1},\dots,i_d} = \sum_{i_j=1}^{I_j} x_{i_1,\dots,i_d} a_{ki_j}, \qquad 1 \leq k \leq m, \qquad j=1,\dots,N.$$
Note that we can also express the mode-$j$ product as the product of two matrices, as $\mat{Y}_{(j)} = \mat{AX}_{(j)}$.

\paragraph{Tucker representation} The Tucker form of an $N$-mode tensor $\ten{X}$, given a target rank $(r_1,\dots,r_N)$, consists of a core tensor $\ten{G} \in \mb{R}^{r_1,\dots, r_N}$ and $N$ factor matrices $\{\mat{A}_j\}_{j=1}^N$, with $\mat{A}_j \in \mb{R}^{I_j \times r_j}$, such that $\ten{X} \approx \ten{G} \bigtimes_{j=1}^N \mat{A}_j$. For short, we denote this form as $\ten{X} \approx [\ten{G};\mat{A}_1, \dots, \mat{A}_N]$. 

A popular algorithm for computing a Tucker representation of a tensor is the Higher-order SVD (HOSVD) \cite{de2000multilinear}. In this algorithm, each mode of a tensor $\ten{X}$ is processed separately and independently. For mode $j$, the factor matrix $\mat{A}_j$ is computed by taking the first $r_j$ left singular vectors of $\mat{X}_{(j)}$, where $(r_1,\dots, r_d)$ is the target rank. Then, the core tensor is formed as $\ten{G} = \ten{X} \bigtimes_{j=1}^d \mat{A}_j^\top$. The tensor compression algorithms we will present in Section~\ref{sec:lowrankalgs} are related to the HOSVD algorithm. This can be computationally expensive for large-scale problems, so we propose new randomized algorithms in Section~\ref{sec:lowrankalgs}, extending some of the methods proposed in~\cite{minster2020randomized}.

\subsection{Randomized Matrix Algorithms}
In this section, we review some of the randomized algorithms (for matrices) that will be used in our later algorithms.  The first of these is the randomized range finder, popularized in \cite{halko2011finding}. Given a matrix $\mat{X} \in \mb{R}^{m \times n}$, this algorithm finds a matrix $\mat{Q}$ that estimates the range of $\mat{X}$, i.e., $\mc{R}(\mat{X}) \approx \mc{R}(\mat{Q})$. This is accomplished by drawing a Gaussian random matrix $\mat{\Omega} \in \mb{R}^{n \times (r+p)}$, where $r$ is the desired target rank, and $p\geq 0$ is an oversampling parameter. %Typically, we will use $p \sim 5$. 
We then form the product $\mat{Y} = \mat{X\Omega}$, which consists of linear combinations of random columns of $\mat{X}$. We next compute a thin QR of the result, $\mat{Y} = \mat{QR}$, to obtain the matrix $\mat{Q}$ such that $\mc{R}(\mat{Q})$ is a good approximation to $\mc{R}(\mat{X})$. Taken together, this process gives the approximation $\mat{X} \approx \mat{QQ}^\top \mat{X}$. See steps 1-3 in Algorithm~\ref{alg:rrid} for the details.

\begin{algorithm}[!ht]
\begin{algorithmic}[1]
    \REQUIRE matrix $\mat{X} \in \mb{R}^{m \times n}$, target rank $r$, oversampling parameter $p$ such that $r+p \leq \min\{m,n\}$
    \ENSURE factor matrix $\mat{F}$, $r+p$ selected indices $\mc{J}$
    \STATE Draw standard random Gaussian matrix $\mat{\Omega} \in \mb{R}^{n \times (r+p)}$
    \STATE Multiply $\mat{Y} \leftarrow \mat{X\Omega}$
    \STATE Compute thin QR $\mat{Y} = \mat{QR}$
    \STATE Perform column pivoted QR on $\mat{Q}^\top$ to obtain permutation matrix $\mat{P} = \bmat{\mat{P}_1 & \mat{P}_2}$, where
    \[ \mat{Q}^\top \bmat{\mat{P}_1 & \mat{P}_2 } = \mat{Z} \bmat{\mat{R}_{11} & \mat{R}_{12}}. \]
    Take $\mc{J}$ to be the index set corresponding to the columns of $\mat{P}_1$ (note that $|\mc{J}| = r+p$).
    \STATE Compute $\mat{F} \leftarrow \mat{Q}(\mat{Q}(\mc{J},:))^{-1}$
\end{algorithmic}
\caption{$[\mat{F},\mc{J}] = $ RRID$(\mat{X},r,p)$}
\label{alg:rrid}
\end{algorithm}

The randomized range finder can be combined with subset selection (specifically pivoted QR, which we call the randomized row interpolatory decomposition (RRID). This algorithm, similar to randomized SVD with row extraction \cite[Algorithm 5.2]{halko2011finding}, produces a low-rank approximation to a matrix $\mat{X}$ that exactly reproduces rows of $\mat{X}$. After the matrix $\mat{Q}$ is computed, we use pivoted QR on $\mat{Q}^\top$ to pick out indices of $r+p$ well-conditioned rows of $\mat{Q}$; denote this by $\mc{J}$. We then form the matrix $\mat{F} = \mat{Q}(\mat{Q}(\mc{J},:))^{-1}$, which gives the low-rank approximation to $\mat{X}$ as $\mat{X} \approx \mat{F}\mat{X}(\mc{J},:)$. In practice, we use the Businger-Golub pivoting algorithm~\cite[Algorithm 5.4.1]{golub2013matrix} but in Theorem~\ref{thm:approx} we use the strong rank-revealing QR factorization~\cite{gu1996efficient}, since among the rank-revealing QR algorithms this algorithm has the best known theoretical results. The details of this algorithm are presented in Algorithm~\ref{alg:rrid}, and will be used later in our tensor compression algorithms.

\section{Tensor-based functional approximation}\label{sec:lowrankalgs}
In this section, we present our tensor-based approach to approximate functions. First, we explain our approach in general, which combines multivariate Chebyshev interpolation and tensor compression techniques to approximate a function $f$. We then propose three tensor compression methods that fit within our framework and the natural application of our approach to kernel approximations. Finally, we provide analysis of the computational cost of our approach as well as analysis of the accuracy.

\subsection{Multivariate functional interpolation}\label{ssec:multi_cheb} Let us denote the Cartesian product of the intervals as 
\[ \mc{D} = [\alpha_1,\beta_1] \times \dots \times [\alpha_N,\beta_N] \subset \mb{R}^N,\]
and let $f: \mc{D} \rightarrow \mb{R}$ be a function that we want to approximate using Chebyshev interpolation. We can expand the function using Chebyshev polynomials as 
\[ f(\xi_1,\dots,\xi_N) \approx \sum_{j_1=1}^n \cdots\sum_{j_N=1}^n f(\eta_{j_1}^{(1)},\dots,\eta_{j_N}^{(N)}) \left( \prod_{k=1}^N S_n^{[\alpha_k,\beta_k]}(\eta_{j_k}^{(k)},\xi_k) \right).\]
where $\eta_{j_k}^{(k)}$ are the Chebyshev points of the first kind in the interval $[\alpha_k,\beta_k]$ for $k=1,\dots,N$ and the interpolating polynomials $S_n$ have been defined in~\eqref{eqn:Sn}. Here we have assumed the same number of Chebyshev points in each dimension for ease of discussion, but this is straightforward to generalize. To match notation with~\eqref{eqn:multi_approx}, we have $c_{j_1,\dots,j_N} = f(\eta_{j_1}^{(1)},\dots,\eta_{j_N}^{(N)})$ and $\phi_{j_k}^{(k)}(\cdot) = S_n^{[\alpha_k,\beta_k]}(\eta_{j_k}^{(k)},\cdot)$.

For a fixed $\boldsymbol{\xi}\in \mb{R}^N$, we can express this in tensor form as 
\[ \hat{f}(\xi_1,\dots,\xi_N) =\ten{M} \bigtimes_{k=1}^N \vec{s}_k(\xi_k),\]
where the entries of $\ten{M} \in \mb{R}^{n\times \dots \times n}$ are $\ten{M}_{j_1,\dots,j_N} = f(\eta_{j_1}^{(1)},\dots,\eta_{j_N}^{(N)})$  and 
\[ \vec{s}_k (\xi_k) = \bmat{ S_n^{[\alpha_k,\beta_k]}(\eta_{1}^{(k)},\xi_k) & \dots & S_n^{[\alpha_k,\beta_k]}(\eta_{n}^{(k)},\xi_k) } \in \mb{R}^{1\times n},  \qquad k=1,\dots,N.\] 
For a new point $\boldsymbol{\xi}$, we need only to compute the factor matrices $\vec{s}_k$.  The core tensor $\ten{M}$ need not be recomputed. However, storing the core tensor can be expensive, especially in high dimensions since the storage cost is $n^N$ entries. Therefore, we propose to work with a compressed representation of the core tensor $\ten{M}$. In~\cite[Section 2.5]{dolgov2021functional}, the authors argue that for certain functions, the degree of the multivariate polynomial required for approximating a function to a specified tolerance can be larger than the multilinear rank (in each dimension). For these functions, this justifies a compression of the tensor which has a smaller memory footprint but is nearly as accurate as the multivariate polynomial interpolant.

Using the Tucker format, we can approximate $\ten{M} \approx \tenh{M} = [\ten{G}; \mat{A}_1,\dots,\mat{A}_N]$ where $\ten{G}^{r_1 \times\dots \times r_N}$ and $\mat{A}_k \in \mb{R}^{n\times r_k}$ for $k=1,\dots,N$. Then, 
\[\hat{f}(\xi_1,\dots,\xi_N) \approx \tenh{M} \bigtimes_{k=1}^N \vec{s}_k(\xi_k) =   \ten{G} \bigtimes_{k=1}^N \widehat{\vec{s}}_k (\xi_k), \]
where $\widehat{\vec{s}}_k (\xi_k) = \vec{s}_k(\xi_k)\mat{A}_k \in \mb{R}^{1\times r_k}$ for $k=1,\dots,N$ with entries 
\[ [\widehat{\vec{s}}_k (\xi_k)]_j =   \sum_{i=1}^n a_{ij}S_n^{[\alpha_k,\beta_k]}(\eta_{i}^{(k)},\xi_k) \qquad j = 1,\dots,r_k.\] 
The compressed Tucker representation $\tenh{M}$ can be obtained using any appropriate tensor compression technique (e.g., HOSVD, sequentially truncated HOSVD). However, computing this low-rank approximation using deterministic techniques is computationally expensive and requires many function evaluations. We propose randomized approaches to computing low-multirank Tucker approximations, in order to both avoid the computational costs and reduce the number of function evaluations, associated with non-randomized approaches.

\subsection{Tensor compression techniques}
We propose three new tensor compression methods in the context of compressing core tensor $\ten{M}$ as defined above. The first uses the Randomized Row Interpolatory Decomposition (RRID), or Algorithm~\ref{alg:rrid}, on each mode unfolding of the tensor to obtain a low-rank approximation for each mode. We provide an analysis for this method as well. The second tensor compression method is similar to the first, but works with a subsampled tensor instead, thus reducing the computational cost. Our third method uses the Kronecker product of Gaussian random matrices instead of a single Gaussian random matrix while compressing $\ten{M}$, which reduces the number of random entries needed. These three methods produce structure-preserving approximations in the Tucker format as described in~\cite{minster2020randomized}; this means that the core has entries from the original tensor $\ten{M}$ and the decomposition  (in exact arithmetic) reproduces certain entries of the tensor exactly.

\subsubsection{Method 1: Randomized Interpolatory Tensor Decomposition}
The first compression algorithm we present
uses a variation of the structure-preserving randomized algorithm that we previously proposed in~\cite[Algorithm 5.1]{minster2020randomized} to compress core tensor $\ten{M}$. In mode $1$, we apply Algorithm~\ref{alg:rrid} with target rank $r_t$ and oversampling parameter $p$ to the mode-1 unfolding $\mat{M}_{(1)}$ to obtain the low-rank approximation
\[ \mat{M}_{(1)} \approx \mat{A}_1 \mat{P}_1^\top\mat{M}_{(1)}. \]
Here $\mat{P}_1 \in \mb{R}^{n\times (r_t+p)}$ has columns from the identity matrix corresponding to the index set $\mc{J}_1$. This process is repeated across each mode to obtain the other factor matrices $\mat{A}_2,\dots,\mat{A}_N$ and the index sets $\mc{J}_2,\dots,\mc{J}_N$.  Finally, we have the low multirank Tucker representation
\[ \ten{M} \approx \tenh{M} = [\ten{G};\mat{A}_1,\dots,\mat{A}_N], \qquad \ten{G} = \ten{M}(\mc{J}_1,\dots,\mc{J}_N). \]
In contrast to~\cite[Algorithm 5.1]{minster2020randomized} in which we used sequential truncation, our low-rank approach can be described as the structure-preserving HOSVD algorithm. To analyze the computational cost of this algorithm, assume that the target tensor rank is  $(r_t,\dots,r_t)$. First, we discuss the cost of tensor compression. We apply the RRID algorithm to each mode unfolding of size $n\times n^{N-1}$. This costs $\mc{O}(r_tn^N)$ flops. The total cost is $\mc{O}(Nr_tn^N)$ flops.

\begin{algorithm}[!ht]
\begin{algorithmic}[1]
\REQUIRE Tensor $\ten{M} \in \mb{R}^{n \times \dots \times n}$ of order $N$, tensor target rank $(r_t,\dots,r_t)$, oversampling parameter $p$ such that $r_t+p \leq n$
\ENSURE approximation $\tenh{M} = [\ten{G}; \mat{A}_1,\dots,\mat{A}_N]$
\FOR {$j=1,\dots,N$}
    \STATE Compute $[\mat{A}_j,\mc{J}_j] =$ RRID$(\mat{M}_{(j)},r_t,p)$ 
\ENDFOR
    \STATE Form core tensor $\ten{G} = \ten{M}(\mc{J}_1,\dots,\mc{J}_N)$  
\end{algorithmic}
\caption{Method 1: Randomized Interpolatory Tensor Decomposition}
\label{alg:var1}
\end{algorithm}

\paragraph{Error Analysis}
We now provide a result on the expected error from compressing a tensor $\ten{M}$ using Algorithm~\ref{alg:var1}. Note that this theorem assumes we are using strong rank-revealing QR (sRRQR) in the RRID algorithm instead of pivoted QR.  This proof is very similar to the proof of \cite[Theorem 5.1]{minster2020randomized}, so we leave the proof to Appendix~\ref{sec:appenda}.

\begin{theorem}\label{thm:approx}
Let $\tenh{M} = [\ten{G}; \mat{A}_1, \dots, \mat{A}_N]$ be the rank-$(r_t, \dots, r_t)$  approximation to $N$-mode tensor $\ten{M} \in \mb{R}^{n \times \dots \times n}$ which is the output of Algorithm~\ref{alg:var1}. We use oversampling parameter $p \geq 2$ such that $\ell = r_t+p < n$, and strong rank-revealing QR factorization~\cite[Algorithm 4]{gu1996efficient} with parameter $f=2$. Then the expected approximation error is 

\begin{equation*}
    \mb{E}_{\{\mat{\Omega}\}_{k=1}^N} \|\ten{M}-\tenh{M}\|_F \leq \sum_{j=1}^N \left( g(n,\ell)\right)^j \left(\left(1+\frac{r_t}{p-1} \right) \sum_{i=r+1}^n \sigma_i^2\left(\mat{M}_{(j)}\right)\right)^{1/2},
\end{equation*}
where $g(n,\ell) = \sqrt{1+4\ell(n-\ell)}$.
\end{theorem}
The interpretation of this theorem is that if the singular values of the mode-unfoldings $\mat{M}_{(j)}$ decay rapidly beyond the target rank, then the low-rank tensor approximation is accurate. See the discussion in~\cite[Section 5.2]{minster2020randomized} for further discussion on the interpretation. 

\subsubsection{Method 2: Randomized Interpolatory Tensor Decomposition with Block Selection}\label{ssec:var2}
Our second compression algorithm aims is a variant of the previous approach, which reduces the number of function evaluations and the computational cost, by working with a subsampled tensor rather than the entire tensor. Recall that we have assumed that the number of Chebyshev grid points in each spatial dimension is the same. Suppose we are given an index set $\mc{I} $ with cardinality $|\mc{I}| = n_b$ representing the subsampled fibers. In the first step of the algorithm, we process mode $1$. We consider the subsampled tensor $\ten{X} = \ten{M}(:,\mc{I},\dots,\mc{I}) \in \mb{R}^{n\times n_b \times \dots \times n_b}$ and apply the RRID algorithm to $\mat{X}_{(1)}$ obtain the matrix $\mat{A}_1$ and the index set $\mc{J}_1$. We then have a low-rank approximation to $\ten{M}$ as 
\[ \mat{M}_{(1)} \approx \mat{A}_1 \mat{P}_1^\top\mat{X}_{(1)},  \qquad  \ten{M} \approx \ten{M}(\mc{J}_1,:,\dots,:) \times_1 \mat{A}_1. \]
Here $\mat{P}_1 \in \mb{R}^{n\times (r_t+p)}$ has columns from the identity matrix corresponding to the index set $\mc{J}_1$. This process is repeated across each mode to obtain the other factor matrices $\mat{A}_2, \dots, \mat{A}_N$ and the index sets $\mc{J}_2,\dots,\mc{J}_N$.  Finally, we have the low-rank Tucker representation
\[ \ten{M} \approx \tenh{M} = [\ten{G};\mat{A},\dots,\mat{A}_N], \qquad \ten{G} = \ten{M}(\mc{J}_1,\dots,\mc{J}_N) . \]
 The summary of all steps for the second variation is described in Algorithm~\ref{alg:var2}.

Compared with Method 1, in each mode we are working with the subsampled version of the mode unfolding rather than the entire mode-unfolding. In other words, by setting $n_b=n$ we can recover Method 1. To subsample, we use the nested property of the Chebyshev nodes~\cite{xu2016chebyshev}. More precisely, let the Chebyshev nodes 
\[ \eta_k^{(n)} = \cos \left(\frac{(2k-1)\pi}{2n}\right) \qquad k=1,\dots,n. \]
Then $\eta_k^{(n)} = \eta_{3k+1}^{(3n)}$; i.e., the $k$-th point $\eta_k^{(n)}$ is the $3k+1$th point $\eta_{3k+1}^{(3n)}$. Therefore, assuming that $n$ is a multiple of $3$, to subsample by 1 level we take $\mc{I} = \{2,5,\dots,n-1\}$. To subsample more aggressively, by $l$ levels, we take the indices to be $\mc{I} = \{ 3^{l-1} + 2, 3^{l-1} + 5,\dots, \}$ and $n_b = |\mc{I}|= n/3^l$.

Now consider the computational cost. At the first step, the unfolded matrix $\mat{X}_{(1)}$ has dimensions $n\times n_b^{N-1}$. Applying RRID for each mode requires $\mc{O}(r_tnn_b^{N-1})$ flops; the total cost for the entire tensor decomposition is therefore $\mc{O}(Nr_tnn_b^{N-1})$.  The number of function evaluations are $nNn_b^{N-1} + (r_t+p)^{N}$; in practice, this is much smaller than $n^N$ since we choose $n_b,r_t \ll n$.

\begin{algorithm}[!ht]
\begin{algorithmic}[1]
    \REQUIRE tensor $\ten{M} \in \mb{R}^{n \times \dots \times n}$ of order $N$,  tensor target rank $(r_t,\dots,r_t)$, index set $\mc{I} \subset \{1,\dots,n\}$ with cardinality $n_b$, oversampling parameter $p$ such that $r_t+p \leq \min\{n,n_b^{N-1}\}$.
    \ENSURE approximation $\tenh{M} = [\ten{G}; \mat{A}_1,\dots,\mat{A}_N]$
    \STATE Set $\mat{P} = \bmat{\mat{I}(:, i_1) & \dots &\mat{I}(:,i_{n_b})} \in \mb{R}^{n\times {n_b}}$ to be the columns of the identity matrix, and $i_1,\dots,i_{n_b} \in \mc{I}$ 
    \FOR {$j=1,\dots,N$}
    \STATE Compute subsampled tensor $\ten{X} = \ten{M} \bigtimes_{k=1, k \neq j }^N  \mat{P}^\top$ 
    \STATE Compute $[\mat{A}_j,\mc{J}_j] =$ RRID$(\mat{X}_{(j)},r_t,p)$
    \ENDFOR
    \STATE Form core tensor $\ten{G} = \ten{M}(\mc{J}_1,\dots,\mc{J}_N)$ 
\end{algorithmic}
\caption{Method 2: Randomized Interpolatory Tensor Decomposition with Block Selection}
\label{alg:var2}
\end{algorithm}

\subsubsection{Method 3: Randomized Kronecker Product }
In the third compression algorithm, we are essentially computing RRID along each mode for the tensor compression step; however, the main difference is that random matrix is generated as the Kronecker product of Gaussian random matrices. More specifically, consider $N-1$ independent standard Gaussian random matrices $\mat\Omega_j \in \mb{R}^{n\times (r_t+p)}$. Here $(r_t,\dots,r_t)$ is the target tensor rank and $p$ is the oversampling parameter. In mode 1, we form the tensor 
\[ \ten{X} = \ten{M} \bigtimes_{j=2}^N \mat\Omega_j^\top \in \mb{R}^{n\times (r_t+p)\times \dots\times (r_t+p)} , \]
and compute the left singular vectors of the mode-1 unfolding $\mat{X}_{(1)}$ to obtain the factor matrix $\mat{A}_1$. We call this the Randomized Kronecker Product approach since using the definition of mode products 
\[ \mat{X}_{(1)} = \mat{M}_{(1)} (\mat\Omega_{N-1} \otimes \dots \otimes \mat\Omega_1) \in \mb{R}^{n\times (r_t+p)^{N-1}}.\] 
To obtain the factor matrices $\mat{A}_2,\dots,\mat{A}_N$, we follow a similar procedure along modes $2$ through $N$. Finally, to obtain the core tensor we compute $\ten{G} = \ten{M}\times_1\mat{A}_1^\top \dots\times_N\mat{A}_N^\top$. The rest of the algorithm resembles Method 1 and is detailed in Algorithm~\ref{alg:var3}. The asymptotic cost of Method 3 is the same as Method 1. However, a potential advantage over Method 1 is the reduction in the number of random entries generated. In Method 1, in using a RRID along each mode, we need to generate $n^{N-1}(r_t+p)$ random numbers, whereas in Method 3, we only need to generate $(N-1)n(r_t+p)$ random numbers. 

\begin{algorithm}[!ht]
\begin{algorithmic}[1]
    \REQUIRE tensor $\ten{M} \in \mb{R}^{n \times \dots \times n}$ of order $N$, tensor target rank $(r_t,\dots,r_t)$,  oversampling parameter $p$ such that $r_t+p \leq n$
    \ENSURE approximation $\tenh{M} = [\ten{G}; \mat{A}_1,\dots,\mat{A}_N]$
    \STATE Draw Gaussian random matrices $\mat\Omega_j \in \mb{R}^{n\times (r_t+p)}$
    \FOR {$j=1,\dots,N$}
        \STATE Form $\ten{X} =  \ten{M} \bigtimes_{k=1, k\neq j}^N \mat\Omega_k^\top$
        \STATE Compute thin QR factorization $\mat{X}_{(j)} = \mat{Q}_j\mat{R}_j$
        \STATE Perform column pivoted QR on $\mat{Q}_j^\top$ to obtain permutation matrix $\bmat{\mat{P}_1 & \mat{P}_2 }$
    \[ \mat{Q}_j^\top \bmat{\mat{P}_1 & \mat{P}_2 } = \mat{Z} \bmat{\mat{R}_{11} & \mat{R}_{12}}. \]
       Take $\mc{J}_j$ to be the index set corresponding to the columns of $\mat{P}_1$ (note that $|\mc{J}_j| = r_t+p$). 
        \STATE Compute $\mat{A}_j = \mat{Q}_j (\mat{Q}_j(\mc{J}_j,:)^{-1}) $
    \ENDFOR
    \STATE Compute core tensor $\ten{G} = \ten{M}(\mc{J}_1,\dots,\mc{J}_N)$.
\end{algorithmic}
\caption{Method 3: Randomized Kronecker Product }
\label{alg:var3}
\end{algorithm}

\subsection{Error analysis} \label{ssec:error}
We now analyze the accuracy of the proposed algorithmic framework. Since there are several approximations (interpolation error, and errors due to tensor compression) involved, it is difficult to give sharp estimates of the error, but we aim to provide some insight into the error analysis. 

 We want to first understand the effect of the tensor compression on the overall approximation error. To this end, we use the triangle inequality to bound
\[ | f(\xi_1,\dots,\xi_N)- \tenh{M} \bigtimes_{j=1}^{N} \vec{s}_j(\xi_j)| \leq |f(\xi_1,\dots,\xi_N)- \ten{M} \bigtimes_{j=1}^{N} \vec{s}_j(\xi_j)| + |( \ten{M}-\tenh{M}) \bigtimes_{j=1}^{N} \vec{s}_j(\xi_j)|.\] 
The first term in the error is due to Chebyshev polynomial interpolation and can be bounded using classical results from multivariate approximation theory. Assuming that the function $f$ satisfies the assumptions in~\cite[Theorem 2.2]{gass2018chebyshev}, then by Corollary 2.3 in that same paper $|f(\xi_1,\dots,\xi_N)- \ten{M} \bigtimes_{j=1}^{N} \vec{s}_j(\xi_j)| \leq C\rho^{-n}$, where $C$ and $\rho$ are positive constants with $\rho > 1$. To bound the second term, we use  the Cauchy-Schwartz inequality to obtain
\[ \begin{aligned} |( \ten{M}-\tenh{M}) \bigtimes_{j=1}^{N} \vec{s}_j(\xi_j)| = & \> | \vec{s}_1(\xi_1) ( \ten{M}-\tenh{M})_{(1)} ([\vec{s}_{N}(\xi_{N})]^\top \otimes \cdots    \otimes [\vec{s}_2(\xi_2)]^\top)| \\ \leq & \>  \|\ten{M} - \tenh{M}\|_F \prod_{t=1}^{N}\|\vec{s}_t(\xi_t)\|_2  .\end{aligned}\] 
We then use the vector inequality $\|\vec{s}_t(\xi_t)\|_2 \leq \sqrt{n}\|\vec{s}_t(\xi_t)\|_\infty$; the terms  $\|\vec{s}_t(\xi_t)\|_\infty$ can also be bounded using results from interpolation; since the terms $S_n$ are the same as the Lagrange basis functions defined at the appropriate interpolation nodes, we can bound
\[ \|\vec{s}_j(\xi_j)\|_\infty\leq \sup_{x \in [\alpha_j,\beta_j]} \sum_{i=1}^n |S_n^{[\alpha_j,\beta_j]}(\eta_{i}^{(j)},x)| =: \Lambda_{j,n},  \qquad j=1,\dots,N,\]
where $\Lambda_{j,n}$ is the Lebesgue constant due to interpolation at the nodes $\{\eta_{i}^{(j)}\}_{i=1}^n$.  For Chebyshev interpolating points, as is the case in our application, it is well known that (see~\cite[Chapter 15]{trefethen2013approximation}) for the Chebyshev points $\Lambda_{t,n} \leq \mc{O}(\log n)$ for $t=1,\dots,N$. This gives $\|\vec{s}_1(\xi_1)\|_2 \leq \mc{O}(\sqrt{n} \log n)$, using which we can derive the error bound
 
\[ |f(\xi_1,\dots,\xi_N)- \tenh{M} \bigtimes_{j=1}^{N} \vec{s}_j(\xi_j)| \leq  C\rho^{-n} +  \mc{O}(n^{N/2}\log^{N} n) \|\ten{M}-\tenh{M}\|_F.\]

The important point to note is that there are two sources of error: the error due to Chebyshev interpolation which decreases geometrically with increasing number of Chebyshev nodes, and the error due to tensor compression is $\|\ten{M}-\tenh{M}\|_F$ which is amplified by the term $\mc{O}(n^{N/2} (\log n)^{N})$. In numerical experiments, we have observed that the error due to the tensor compression was sufficiently small, so the amplification factor did not appear to have deleterious effects. 

It is desirable to have an adaptive approach that produces an approximation for a given accuracy; we briefly indicate a few ideas to develop an adaptive approach. However, we do not pursue the adaptive case here further and leave it for future work. If the number of Chebyshev points $n$ is not known in advance, an adaptive approach can be derived that combines the methods proposed here with the computational phases as described in~\cite{hashemi2017chebfun,dolgov2021functional}. We can start with a small estimate for $n$, we can use the fact that the Chebyshev points of the first kind are nested. We can start with a small value of $n$ (say $n= 16$), and increase $n$ by multiples of $3$ until a desired accuracy is achieved (by evaluating the accuracy on a small number of random points); since the points are nested, this avoids recomputing a fraction of the function evaluations. If the target rank $r_t$ is not known in advance, we can use the techniques described in~\cite[Section 4]{minster2020randomized} to estimate the target rank such that the low-rank decomposition satisfies a relative error bound involving a user-specified tolerance.

\section{Application: Low-rank Kernel approximations}\label{sec:kernel}
 In this section, we use the functional tensor approximations developed in Section~\ref{sec:lowrankalgs} to compute low-rank approximations to kernel matrices.
\subsection{Problem Setup and Notation}\label{ssec:probsetup}
Rank-structured matrices are a class of dense matrices that are useful for efficiently representing kernel matrices. In this approach, the kernel matrix is represented as a hierarchy of subblocks, and certain off-diagonal subblocks are approximated in a low-rank format. There are several different types of hierarchical formats for rank-structured matrices, including $\mc{H}$-matrices, $\mc{H}^2$-matrices, hierarchically semi-separable (HSS) matrices, and hierarchically off-diagonal low-rank (HODLR) matrices. For more details on hierarchical matrices, see \cite{hackbusch1999sparse,hackbusch2002data,grasedyck2003construction,borm2003hierarchical,hackbusch2015hierarchical,sauter2010boundary}. If the number of interaction points is $N_s$, then the matrix can be stored approximately using $\mc{O}(N_s(\log^s N_s))$ entries rather than $N_s^2$ entries, and the cost of forming matrix-vector products (matvecs) is $\mc{O}(N_s(\log^t N_s))$ flops, where $s,t\geq 0$ are nonnegative integers depending on the format and the algorithm for computing the rank-structured matrix. The dominant cost in working with rank-structured matrices is the high cost of computing low-rank approximations corresponding to off-diagonal blocks, which we address here using novel tensor-based methods.

\begin{figure}[!ht]
    \centering
    \includegraphics[scale=.4]{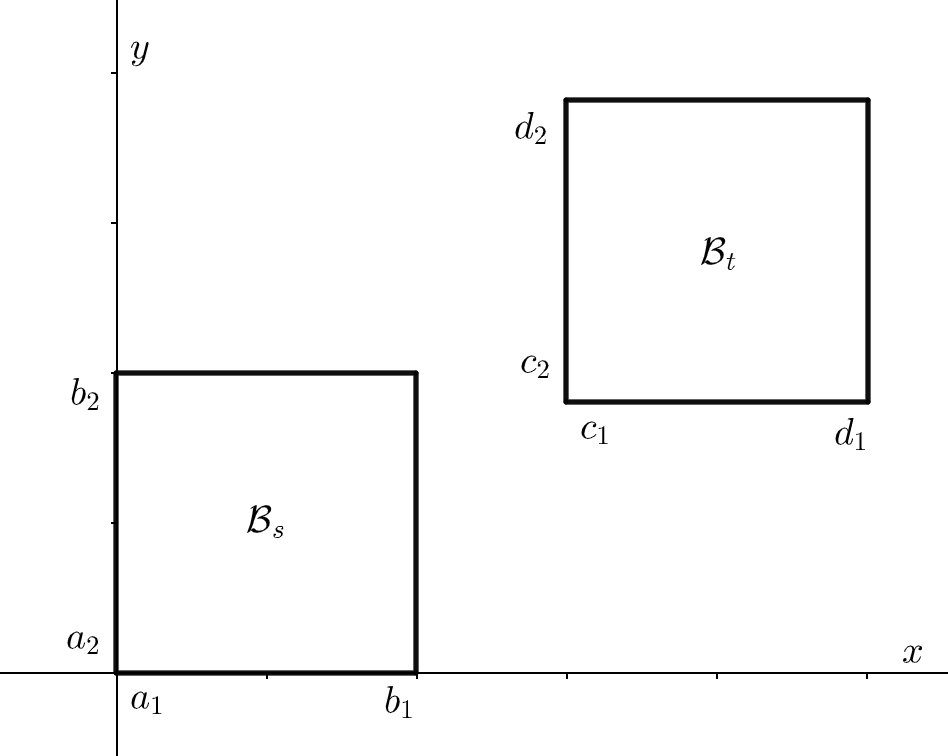}
    \caption{Visual of the bounding boxes in two spatial dimensions for source points, $\mc{B}_s$, and for target points, $\mc{B}_t$. Note that $\mc{X} \subset \mc{B}_s$ and $\mc{Y} \subset \mc{B}_t$.}
    \label{fig:bounds}
\end{figure}

Before showing how to compute these off-diagonal blocks efficiently, we define the problem setup that we will use in the rest of this chapter. Let $\mc{X} = \{\mat{x}_1,\mat{x}_2,\dots,\mat{x}_{N_s}\}$ and  $\mathcal{Y} = \{\mat{y}_1,\mat{y}_2,\dots,\mat{y}_{N_t}\}$ be $N_s$ and $N_t$ source and target points, respectively. These points are assumed to be enclosed by two boxes, $\mc{B}_s, \mc{B}_t \subset \mb{R}^D$, where $D$ is the dimension, defined as
\begin{equation}\label{eq:blockdef}
    \begin{aligned}
    \mc{B}_s &= [a_1,b_1]\times \dots \times [a_D,b_D], \\
    \mc{B}_t &= [c_1,d_1] \times \dots \times[c_D,d_D],
    \end{aligned}
\end{equation}
such that $\mc{X} \subset \mc{B}_s$ and $\mc{Y} \subset \mc{B}_t$. These bounding boxes are plotted in Figure~\ref{fig:bounds} for visualization in two spatial dimensions (i.e., $D=2$).
Then, for the kernel $\kappa: \mb{R}^D \times \mb{R}^D \to \mb{R}$, let the interaction matrix $\matc{K}(\mathcal{X},\mathcal{Y})$ be defined as
\begin{equation}\label{eq:kmat}
\matc{K}(\mc{X},\mc{Y}) = \bmat{\kappa(\mat{x}_1,\mat{y}_1) & \kappa(\mat{x}_1,\mat{y}_2) & \dots & \kappa(\mat{x}_1,\mat{y}_{n_t}) \\ \kappa(\mat{x}_2,\mat{y}_1) & \kappa(\mat{x}_2,\mat{y}_2) & \dots & \kappa(\mat{x}_2, \mat{y}_{N_t}) \\ \vdots & \vdots & \ddots & \vdots \\ \kappa(\mat{x}_{N_s},\mat{y}_1) & \kappa(\mat{x}_{N_s},\mat{y}_2) & \dots & \kappa(\mat{x}_{N_s},\mat{y}_{N_t})}.
\end{equation}

We can only represent blocks in low-rank form if those blocks are \textit{admissible}.  
Note that the diameters of the bounding blocks $\mc{B}_s$ and $\mc{B}_t$ are
\begin{equation}\label{eq:diam}\text{diam}(\mc{B}_s) = \sqrt{\sum_{j=1}^D(b_j-a_j)^2}, \qquad \text{diam}(\mc{B}_t) = \sqrt{\sum_{j=1}^D(d_j-c_j)^2},
\end{equation}
We define the distance $\text{dist}(\mc{B}_s,\mc{B}_t)$ between these blocks as the minimum distance between any two points from $\mc{B}_s$ and $\mc{B}_t$. Specifically,
\begin{equation}\label{eq:blockdist}
    \text{dist}(\mc{B}_s,\mc{B}_t) = \min_{\substack{\mat{x}\in \mc{B}_s \\ \mat{y} \in \mc{B}_t}} \left( \sum_{j=1}^D(y_j-x_j)^2 \right)^{1/2}.
\end{equation}
The domains $\mc{B}_s$ and $\mc{B}_t$ are considered \textit{strongly admissible} if they satisfy 
\begin{equation}\label{eq:admissible}
 \max\{\text{diam}(\mc{B}_s),\text{diam}(\mc{B}_t)\} \leq \eta \, \text{dist}(\mc{B}_s,\mc{B}_t), 
\end{equation}
where $\eta > 0$ is the admissibility parameter. This means two sets of points are strongly admissible if they are sufficiently well-separated. Two blocks are \textit{weakly admissible} if there is no overlap between them. Typically, we will assume that the blocks $\mc{B}_s$ and $\mc{B}_t$ are either strongly or weakly admissible. However, in Section~\ref{ssec:gauss} we discuss how to compute a global low-rank approximation assuming only that $B_s = B_t$ and that the kernel is sufficiently smooth (and that the blocks are not admissible).

Given points $\mc{X} \subset \mc{B}_s$ and $\mc{Y} \subset \mc{B}_t$, we seek to compute the pairwise kernel interactions between the sets of points, i.e., $\matc{K}(\mc{X},\mc{Y})$. This is expensive both from a computational and storage perspective. If we simply compute the kernel interaction for each pair of points, this requires $N_s N_t$ kernel interactions and the cost of storage is $\mc{O}(N_sN_t)$.  As $N_s$ and $N_t$ can get quite large, computing $\matc{K}(\mc{X},\mc{Y})$ can become computationally challenging. Assuming that the points are well-separated so that the kernel interactions between the source and targets can be treated as a smooth function, the matrix $\matc{K}(\mc{X},\mc{Y})$ can be approximated in a low-rank format. However, computing this low-rank approximation is computationally expensive. To this end, our goal is to compute an approximation to  $\matc{K}(\mc{X},\mc{Y})$ with a cost that is linear in $N_s$ and $N_t$, which we accomplish using multivariate Chebyshev interpolation.

\begin{figure}[!ht]
    \centering
    \includegraphics[scale=.4]{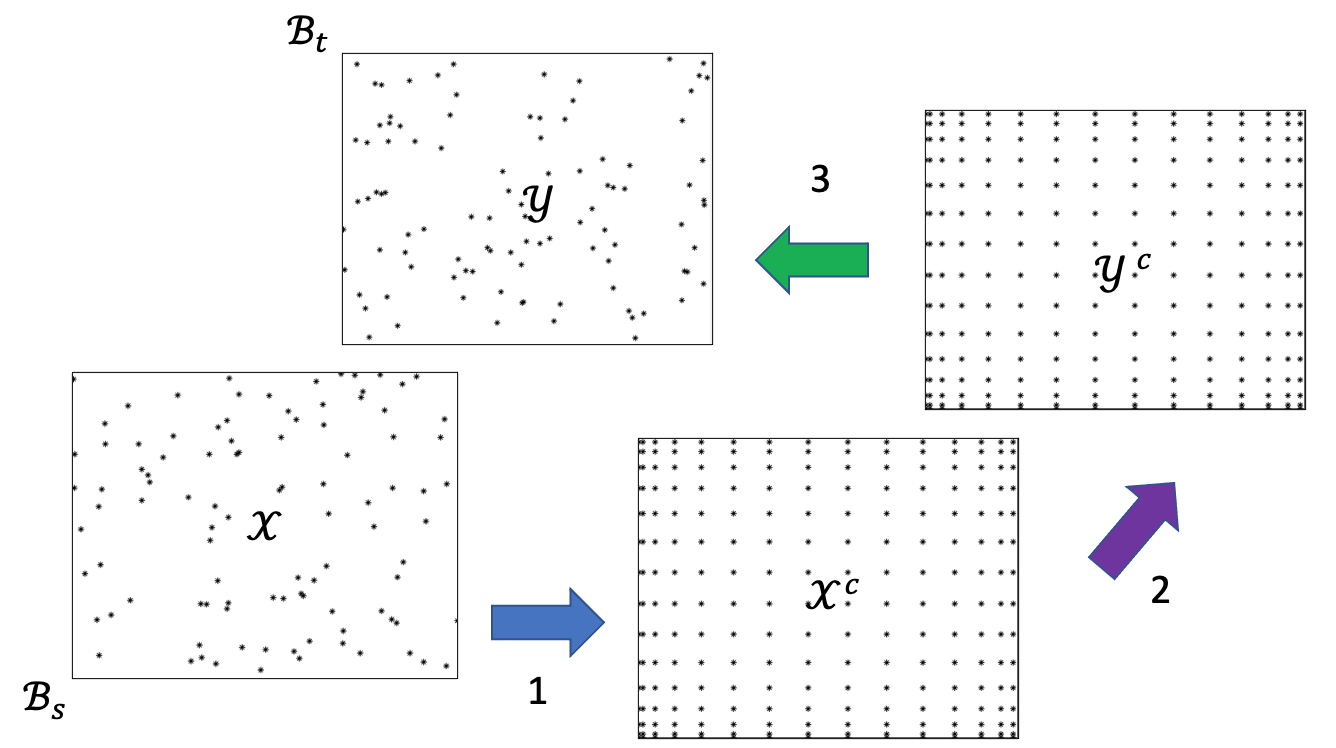}
    \caption{Visual representation in two spatial dimensions of mapping source and target points to Chebyshev grids, shown by steps 1 and 3. Step 2 shows computing the interactions between Chebyshev grid. Each numbered step corresponds to a matrix shown in Figure~\ref{fig:rankmats}. }
    \label{fig:cheby}
\end{figure}

\subsection{Low-rank approximation of $\matc{K}(\mc{X},\mc{Y})$}\label{ssec:lowrank}
Suppose we have the boxes $\mc{B}_s,\mc{B}_t \subset \mb{R}^{D}$ denoting the bounding boxes for the sources and the targets respectively that are well-separated so that the kernel function $\kappa: \mc{B}_s \times \mc{B}_t \rightarrow \mb{R}$ can be viewed as a smooth function $f:\mb{R}^{2D} \rightarrow \mb{R}$ of the form 
\[ f(x_1,\dots,x_D, y_1,\dots,y_D) := \kappa(\vec{x},\vec{y}), \qquad \vec{x} \in \mc{B}_s, \vec{y}\in \mc{B}_t.\]
That is, we identify $\xi_j = x_j$ and $\xi_{j+D} = y_j$ for $j=1,\dots,D$. Similarly, $\alpha_j = a_j, \beta_j = b_j, \alpha_{j+D} = c_j$ and $\beta_{j+D} = d_j$ for $j=1,\dots,D$. This allows us to approximate the function $f$, and, therefore, the kernel by a multivariate Chebyshev polynomial approximation
\[ \begin{aligned} \kappa(\vec{x},\vec{y}) \approx & \> \sum_{j_1=1}^n \cdots\sum_{j_N=1}^n f(\eta_{j_1}^{(1)},\dots,\eta_{j_N}^{(N)}) \left( \prod_{k=1}^N S_n^{[\alpha_k,\beta_k]}(\eta_{j_k}^{(k)},\xi_k) \right) \\
= & \> \ten{M} \left( \bigtimes_{k=1}^D  \mat{s}_k(x_k) \right)  \left(\bigtimes_{k=D+1}^{2D} \mat{s}_k(y_{k-D})\right), 
\end{aligned}
\]
following the notation in Subsection~\ref{ssec:multi_cheb} for $\ten{M}$. We can use the multidimensional Chebyshev approximation to obtain a low-rank representation for $\matc{K}(\mc{X},\mc{Y})$ as follows.

\begin{figure}[!ht]
    \centering
    \includegraphics[scale=.4]{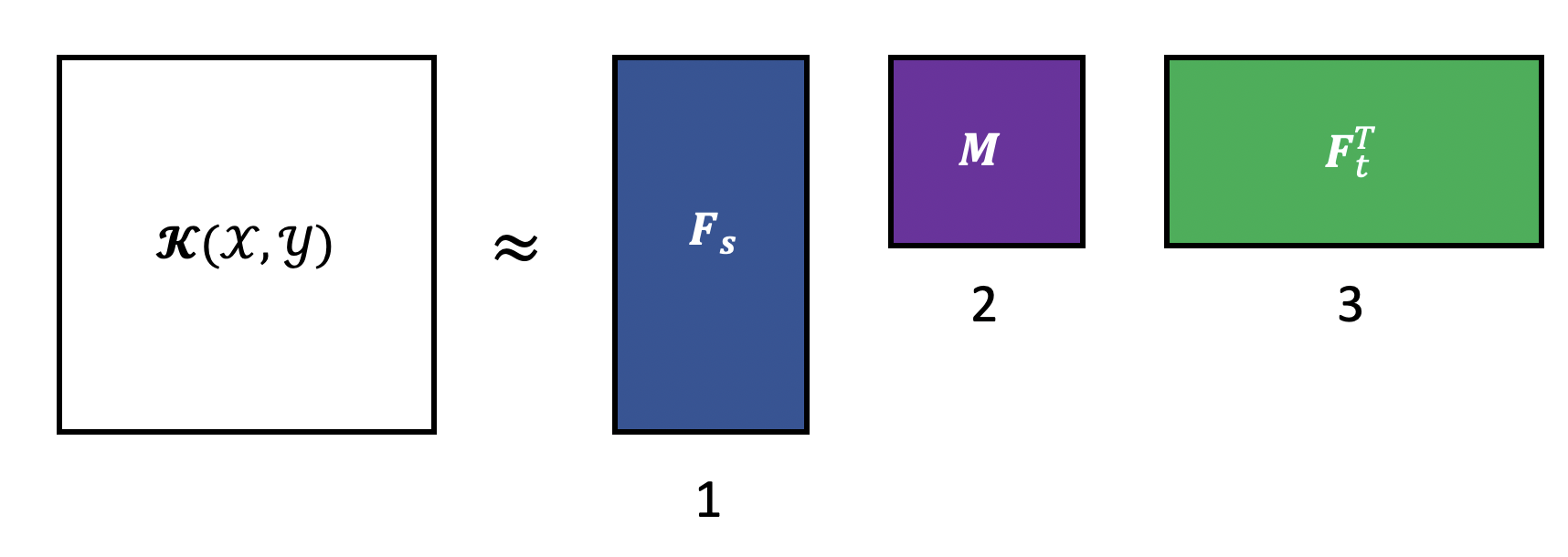}
    \caption{Visualization of the matrices comprising the interaction matrix approximation process, corresponding to steps shown in Figure~\ref{fig:cheby}. The matrix $\mat{F}_s$ maps the source points to a Chebyshev grid, matrix $\mat{F}_t$ maps the target points to a Chebyshev grid, and $\mat{M}$ computes the interactions between Chebyshev grids.}
    \label{fig:rankmats}
\end{figure}

\paragraph{Low-rank approximation} We are now given the set of points $\mc{X} = \{ \vec{x}_j\}_{j=1}^{N_s}$ and $\mc{Y} = \{ \vec{y}_j\}_{j=1}^{N_t}$. Let us define the sequence of matrices $\{\mat{U}_j\}_{j=1}^D$ and $\{\mat{V}_j\}_{j=1}^D$ such that 
\[\mat{U}_j = \> \bmat{\vec{s}_j([\vec{x}_1]_j) \\ \vec{s}_j([\vec{x}_2]_j)] \\ \vdots \\ \vec{s}_j([\vec{x}_{N_s}]_j)} \in \mb{R}^{N_s \times n}, \qquad \mat{V}_j =  \bmat{\vec{s}_{j+D}([\vec{y}_1]_j) \\  \vec{s}_{j+D}([\vec{y}_2]_j) \\ \vdots \\ \vec{s}_{j+D}([\vec{y}_{N_t}]_j)}^\top \in \mb{R}^{N_t \times n},\]
for $j=1,\dots,D$.  Let $\ltimes$ define the row-wise Khatri-Rao product. Using these sequences, we define the following factor matrices to be the row-wise Khatri-Rao product of the matrices $\{\mat{U}_j\}_{j=1}^D$ and $\{\mat{V}_j\}_{j=1}^D$ in reverse order
\[ \begin{aligned}
    \mat{F}_s = &\> \mat{U}_{D} \ltimes \mat{U}_{D-1} \ltimes \dots \ltimes \mat{U}_1 \in \mb{R}^{N_s \times n^D} \\
    \mat{F}_t = & \> \mat{V}_{D} \ltimes \mat{V}_{D-1} \ltimes \dots \ltimes \mat{V}_{1} \in \mb{R}^{N_t \times n^D}.
\end{aligned}\]
This gives the low-rank approximation (see Figure~\ref{fig:rankmats} for a visual representation)
\begin{equation}\label{eqn:lowrank1}
    \matc{K}(\mc{X},\mc{Y}) \approx \mat{F}_s \mat{M} \mat{F}_t^\top, 
\end{equation} 
where $\mat{M} = \mat{M}_{(1:D)} \in \mb{R}^{n^D\times n^D}$ is an unfolding of the tensor $\ten{M}$. The storage cost of this representation is $n^{2D} + nD(N_s+N_t)$, assuming $\mat{F}_s$ and $\mat{F}_t$ are stored in terms of factors  $\{\mat{U}_j\}_{j=1}^D$, $\{\mat{V}_j\}_{j=1}^D$ and are not formed explicitly.

\paragraph{Compressed low-rank representation} Suppose we compress $\ten{M}$ using the methods discussed in Section~\ref{sec:lowrankalgs} to obtain $\tenh{M} = [\ten{G}; \mat{A}_1, \dots, \mat{A}_{2D}]$ with multirank $(\ell,\dots,\ell)$ (here $\ell = r_t + p$), we can write the kernel approximation as 
$$\kappa(\mat{x},\mat{y}) \approx \tenh{M} \bigtimes_{k=1}^{2D} \mat{s}_k(\xi_k) = \ten{G} \bigtimes_{k=1}^{2D} \mathat{s}_k(\xi_k),$$
where $\mathat{s}_k(\xi_k) = \mat{s}_k(\xi_k) \mat{A}_k$ for $k = 1,\dots, 2D$. Similar to the low-rank representation discussed above, we define the sequences $\{\mathat{U}_j\}_{j=1}^D$, $\{\mathat{V}_j\}_{j=1}^D$, as 
\[ \mathat{U}_j = \mat{U}_j\mat{A}_j \qquad \mathat{V}_j =  \mat{V}_j\mat{A}_j \qquad j=1,\dots,D.\]
We also have the matrices
\[ \begin{aligned}
    \mathat{F}_s = &\> \mathat{U}_{D} \ltimes \mathat{U}_{D-1} \ltimes \dots \ltimes \mathat{U}_1 \in \mb{R}^{N_s \times \ell^D} \\
    \mathat{F}_t = & \> \mathat{V}_{2D} \ltimes \mathat{V}_{2D-1} \ltimes \dots \ltimes \mathat{V}_{D+1} \in \mb{R}^{N_t \times \ell^D}.
\end{aligned}\]
Finally, the low-rank matrix approximation to $\matc{K}(\mc{X},\mc{Y})$ is 
\begin{equation}\label{eqn:lowrank2}
    \matc{K}(\mc{X},\mc{Y}) \approx \mathat{F}_s \mathat{M} \mathat{F}_t^\top,
\end{equation}
where, now, $\mathat{M} = \mat{G}_{(1:D)} \in \mb{R}^{\ell^D \times \ell^D}$. The storage cost using this representation (assuming $\mat{F}_s$ and $\mat{F}_t$ are stored in terms of factors  $\{\mathat{U}_j\}_{j=1}^D$, $\{\mathat{V}_j\}_{j=1}^D$) is $\ell^{2D} + \ell D(N_s+N_t)$ and has obvious storage benefits when $\ell\ll n$. 

\paragraph{Further compression} We now have a low-rank approximation of the form $\mathat{F}_s \mathat{M} \mathat{F}_t^\top$. For an additional computational cost, we can achieve further compression. We compute the thin QR factorizations $\mathat{F}_s = \mat{Q}_s\mat{R}_s$ and $\mathat{F}_t = \mat{Q}_t\mat{R}_t$ and form $\mat{B} = \mat{R}_s\mathat{M} \mat{R}_t^\top$.  Then we compute the truncated SVD $\mat{B} \approx \mathat{U}_B \mathat{\Sigma}_B \mathat{V}_B^\top$ corresponding to the matrix rank $r_t$. Finally, we have the low-rank approximation 
\begin{equation}\label{eqn:lowrank3}\matc{K}(\mc{X},\mc{Y}) \approx \mathat{U}\mathat{\Sigma}\mathat{V}^\top, 
\end{equation}
where $\mathat{U}  = \mat{Q}_s \mathat{U}_B$ and $\mathat{V}  = \mat{Q}_t\mathat{V}_B$. Suppose the target matrix rank is $r_m$. Assuming that we use RandSVD to compute the SVD of $\mat{B}$, the additional cost required for further compression is $\mc{O}((r_mr_t + r_t^{2D})(N_s+N_t) + r_t^{2D}r_m)$ flops.

\subsection{Computational cost and Error}
Suppose we were to use the SVD to compress interaction matrix $\matc{K}(\mc{X},\mc{Y})$; the cost of compression is then $\mc{O}(\max\{N_s,N_t\}\min\{N_s,N_t\}^2)$ flops. If the number of sources and targets are large, then this cost has cubic scaling with the number of points. By using the BBFMM, the cost of forming the kernel approximation $\matc{K}(\mc{X},\mc{Y}) \approx \mat{F}_S\mat{MF}_T^\top$ becomes $\mc{O}(n^D(N_s+N_t) + n^{2D})$ flops, since forming $\mat{F}_S,\mat{F}_T$ cost $\mc{O}(n^D(N_s+N_t))$ flops and forming $\mat{M}$ costs $\mc{O}(n^{2D})$ flops. Since the cost of storage of the BBFMM approximation can be large when $n$ is large, we have proposed various techniques for compressing $\mat{M}$. 

We now summarize the costs of the various compression techniques at our disposal. It should be noted that to obtain a low-rank approximation to $\matc{K}(\mc{X},\mc{Y})$, we need to include the cost of forming $\mat{F}_s$ and $\mat{F}_t$ which cost $\mc{O}(n^D(N_s+N_t))$ flops. If the SVD is used for compressing $\mat{M}$, then the cost is $\mc{O}(n^{3D})$ flops; if RandSVD is used instead of SVD, then the cost is $\mc{O}(n^{2D}r_m )$ flops. Methods 1 and 3 both cost $\mc{O}\left(Dr_tn^{2D} \right)$ flops, whereas Method 2 costs $\mc{O}(Dn_b^{2D-1}n r_t)$ flops. Method 2 is the most computationally efficient of the proposed algorithms if $n_b,r_t,p \ll n$, and if $n_b=n$, Method 2 has the same cost as Methods 1 and 3. Note that Methods 1 and 3 have comparable cost to RandSVD directly computed on $\mat{M}$, but typically the tensor rank $r_t$ is chosen to be smaller than $r_m$, so we expect the tensor-based methods to be more efficient for the same accuracy. As mentioned earlier, Method 3 may be preferable to Method 1 since it requires generating fewer random numbers. 

\begin{table}[!ht]
    \centering
    \begin{tabular}{c|c|c|c}
         Approximation Method & Cost & Kernel Evals. & Storage \\
         \hline
         SVD  & $\mc{O}(\max\{N_s,N_t\}\min\{N_s,N_t\}^2)$ & $N_sN_t$ & $\mc{O}(r(N_s+N_t))$\\
         BBFMM & $\mc{O}(n^D(N_s+N_t) + n^{2D})$ & $n^{2D}$ & $\mc{O}(n^{2D} + n^D N_sN_t)$ \\
        BBFMM + SVD & $\mc{O}(n^{2D} + n^D(N_s+N_t))$ & $n^{2D}$ & $\mc{O}(r_m(N_s+N_t))$
    \end{tabular}
    \vspace{.3cm}
    
    \begin{tabular}{c|c|c|c}
        Method & Cost &  Kernel Evals. & Storage\\
        \hline
        RandSVD & $\mc{O}(r_mn^{2D} + r_mn^D(N_s+N_t) )$ & $n^{2D}$ & $\mc{O}(r_m(N_s+N_t))$\\
         Method 1 & $\mc{O}\left(Dr_tn^{2D} + (r_t^D + Dr_tn)(N_s+N_t)\right)$ & $n^{2D}$ & $\mc{O}(r_t^{2D} + r_t^D(N_s+N_t))$\\
        Method 2 & $\mc{O}(Dn_b^{2D-1}n r_t + (r_t^D + Dr_tn)(N_s+N_t))$ & $2Dn_b^{2D-1}n +\ell^{2D}$ & $\mc{O}(r_t^{2D} + r_t^D(N_s+N_t))$\\
        Method 3 & $\mc{O}\left(Dr_tn^{2D} + (r_t^D + Dr_tn)(N_s+N_t)\right)$ & $n^{2D}$ & $\mc{O}(r_t^{2D} + r_t^D(N_s+N_t))$ \\
    \end{tabular}
    \caption{Summary of the computational cost of the various algorithms. The upper table represents the computational cost of the standard approaches, and the lower table represents the computational costs associated with the tensor compression methods proposed. Here, $n$ is the number of Chebyshev nodes, $r_t$ is the tensor target rank, $r_m$ is the matrix target rank, $\ell = r_t + p$ where $p$ is the oversampling parameter, and $n_b$ is the number of blocks drawn from $\ten{M}$ in Method 2.}
    \label{tab:compcost}
\end{table}

Another potential benefit of Method 2 is the number of kernel evaluations needed to form $\ten{M}$. In nearly all the other approaches, the tensor $\ten{M}$ needs to be formed explicitly requiring $n^{2D}$ kernel evaluations. In Method 2, only a portion of the tensor $\ten{M}$ needs to be formed. That is, for $D=2$ (i.e., in two spatial dimensions),  we only need  $4n_b^3n +(r_t+p)^4$ kernel evaluations. Assuming $r_t = p = 5$ and $n = 27$, we only need $(4n_b^3n +(r_t+p)^4)/n^4 \approx 15\% $ of the total function evaluations if $n_b = 9$ and approximately $2\%$ of the total function evaluations if $n_b = 3$. In other words, we only need a small percentage of the kernel evaluations, which can be beneficial if the kernel evaluations are expensive.

There are other benefits to using our tensor-based compression method. Consider a time-dependent problem in which the source and target points are changing in time but are contained within the same bounding boxes $\mc{B}_s$ and $\mc{B}_t$. Although the matrices $\mat{F}_s,\mat{F}_t$ change in time, the kernel interactions in $\ten{M}$ do not change, and the matrix approximation $\tenh{M}$ can be precomputed and deployed efficiently in the time dependent problem. This can yield significant reductions in computational cost and has low storage costs since the resulting approximations can be stored efficiently.

\paragraph{Error analysis} 

Consider a fixed source point $\vec{x} \in \mc{B}_s$ and a fixed target point $\vec{y}\in \mc{B}_t$. With the notation in Section~\ref{ssec:lowrank}, we can write $\kappa(\vec{x},\vec{y}) \approx \ten{M} \bigtimes_{j=1}^{2D} \vec{s}_j(\xi_j)$, where $$\vec{s}_j(\xi_j) = \bmat{S_n^{[a_j,b_j]}(\eta_{1}^{(j)},\xi_j) &\cdots&  S_n^{[a_j,b_j]}(\eta_{n}^{(j)},\xi_j)} \in \mb{R}^{1\times n}, \qquad j = 1,\dots,2D. $$ 
Note that $\xi_j = x_j$  and $\xi_{j+D} = y_j$ for $j=1,\dots,D$. 
The vectors $\vec{s}_{j+D}(\xi_{j+D})$ for $j=1,\dots,D$ are defined analogously. Following the discussion in Section~\ref{ssec:error}, we can derive the error bound
 
\[ |\kappa(\vec{x},\vec{y})- \tenh{M} \times_{j=1}^{2D} \vec{s}_j(\xi_j)| \leq  C\rho^{-n} +  \mc{O}(n^D\log^{2D} n) \|\ten{M}-\tenh{M}\|_F.\]

\section{Numerical Experiments}\label{sec:results}
We demonstrate the performance of our algorithms on three different sets of numerical experiments. The first set of experiments study the accuracy of the tensor-based function approximations on synthetic test problems and an application in electrical circuits. The second set of experiments study the performance of the low-rank approximations for kernel matrices on a wide variety of kernels from PDEs, Radial basis interpolation, and Gaussian processes. The final set of experiments, is an application to a test dataset that has been used in Gaussian processes.

\subsection{Function approximation}
We first investigate the use of the tensor-based function approximation. We choose three different test problems that were also considered in~\cite{hashemi2017chebfun,dolgov2021functional}. Note that while these references focus on constructing highly accurate function approximations (typically with $15$ digits of accuracy), our goal is to test the performance of the tensor compression methods.  Specifically, we choose
\begin{equation}\label{eqn:testfns}
    \begin{aligned}
    f_1(x,y,z) = & \> \frac{1}{1+25(x^2 + y^2 + z^2)} \\ 
    f_2(x,y,z) = & \> \sin(x+yz)\\
    f_3(x,y,z) = & \> \text{tanh}(3(x+y+z)).
    \end{aligned}
\end{equation}
We choose the number of Chebyshev nodes $n=36$, the target rank $\ell = r_t + p = 10$. We choose the number of blocks $n_b = 4$. To compute the error, we used $100$ randomly generated points in $[-1,1]^3$ and report the relative error in the $\infty$-norm. The results are reported in Table~\ref{tab:synth}. We see that Method 3 is the most accurate and its accuracy is comparable with Method 1 and HOSVD. For functions $f_1$ and $f_3$, the low-rank compression is not that significant compared to function $f_2$ and for these cases the performance of Method 2 are comparable to Method 3 (and HOSVD). For $f_2$, Method 2 has small error but is not as accurate as the other methods. However, Method 2 requires far fewer function evaluations $\approx 5 \%$ and can be computationally efficient for good accuracy when the function evaluations are very expensive. 
\begin{table}[!ht]
    \centering
    \begin{tabular}{c|c|c|c|c}
        &HOSVD &  Method 1 & Method 2 & Method 3  \\ \hline 
    $f_1$     & $8.75\times 10^{-3}$& $8.75\times 10^{-3} $ & $8.75\times 10^{-3} $ & $2.29\times 10^{-3} $ \\ 
    $f_2$ &  $6.49\times 10^{-13}$ & $5.80\times 10^{-12}$ & $1.046\times 10^{-8}$ & $2.41\times 10^{-13}$ \\ 
    $f_3 $ & $2.71\times 10^{-3}$ & $7.18\times 10^{-2}$ & $4.41\times 10^{-2}$ & $5.00\times 10^{-3}$
    \end{tabular}
    \caption{Accuracy of the proposed methods on the test functions defined in~\eqref{eqn:testfns}.}
    \label{tab:synth}
\end{table}

\paragraph{A six-dimensional test problem} Our next test problem is a 6-dimensional test problem modeling an electrical circuit, which has been used for statistical screening and sensitivity analysis~\cite{constantine2017global}. The details of this test problem are available in~\cite{sfuwebsite}, under Emulation \& Prediction test problems/OTL Circuit test problem. For the test parameters, we choose the number of Chebyshev nodes $n=12$, the rank $r_t = 5$ and oversampling $p=0$. We choose the number of blocks $n_b = 4$.  To compute the error, we used $100$ randomly generated points in the domain of the function, and report the relative error in the $\infty$-norm:
\begin{center}
\begin{tabular}{c|c|c|c}
        HOSVD &  Method 1 & Method 2 & Method 3  \\ \hline
        $7.74\times 10^{-8}$ & $2.04\times 10^{-7}$ & $1.77\times 10^{-7}$ &$1.83\times 10^{-7}$.
        \end{tabular}
\end{center}
For this example, all three methods had comparable accuracy and Method 2 only requires $\sim 3\%$ of the function evaluations. The resulting functional approximation using any of the methods can, therefore, be deployed as a numerical surrogate in applications. Method 2 is especially attractive for problems in which the function evaluations are expensive.  

\subsection{Kernel low-rank approximations}

First, we consider the kernel approximations in 2D. In these experiments, we use $n = 27$ Chebyshev nodes for each dimension.   For Method 2, we will take $n_b = 9$ for all the kernels. Also, we will compute the relative error as 
\[ \text{relerr} = \frac{\|\matc{K}(\mc{X},\mc{Y}) - \mathat{F}_s\mathat{M}\mathat{F}_t^\top\|_{\max}}{\|\matc{K}(\mc{X},\mc{Y})\|_{\max}},\]
where $\|\mat{A}\|_{\max} = \max_{i,j}|a_{ij}|$ is the max-norm. 

\begin{figure}[!ht]
    \centering
    \includegraphics[scale=.4]{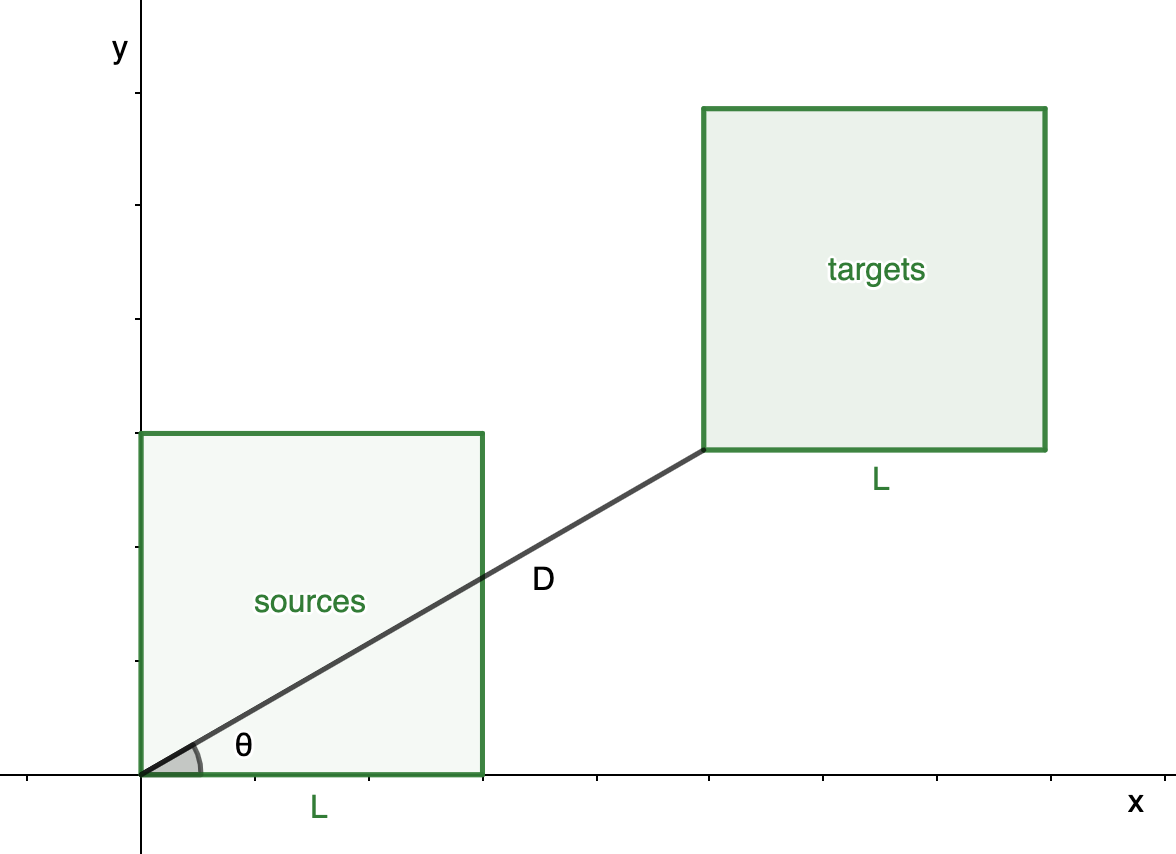}
    \caption{Box setup for our numerical experiments for two spatial dimensions, where $L$ is the length of both the source and target boxes, $D$ is the distance between bottom left corners of boxes, and $\theta$ is the angle describing the placement of the target box.}
    \label{fig:boxsetup}
\end{figure}

We define our source and target boxes for two spatial dimensions in the following way. Let the source box have one vertex at the origin and side length $L$. The target box also has side length $L$, and let $D$ be the distance between bottom left vertices. Finally, define $\theta$ to be the angle describing the placement of the targets box in relation to the sources box. A visual representation of this setup is in Figure~\ref{fig:boxsetup}. Unless otherwise stated, we will use $L = 5$, $D = 10$, and $\theta = \pi/4$. In the boxes defined in this manner, we will generate $N_s = 500$ and $N_t = 500$ uniformly randomly distributed source and target points, respectively.

\paragraph{Experiment 1: Comparing different kernels} In our first experiment, we consider the relative error produced by our algorithms for five different kernels that are taken from different fields: partial differential equations (PDEs), Radial basis functions, and Gaussian Processes.
\begin{table}[!ht]\centering
\begin{tabular}{c|c|c}
    Kernel function $\varphi(r)$ & Name & Field \\ \hline
   $1/r$ & Laplace-3D &  PDEs \\
    $1/r^2$ & Biharmonic &  PDEs \\
    $-\log(r)$ & Laplace-2D&  PDEs \\
    $r^2\log(r)$ & Thin-plate&  Radial basis functions \\
       $\sqrt{1+(\frac{r}{\sigma})^2}$ & Multiquadric &   Radial basis functions \\
    $\exp(-(\frac{r}{\sigma})^2)$ & Gaussian &   Radial basis functions \\
       $\exp(-\frac{r}{\sigma})$ & Mat\'ern-$1/2$ &  Gaussian Processes \\
    $(1+\sqrt{3}(\frac{r}{\sigma}))\exp(-\frac{\sqrt{3}r}{\sigma})$ & Mat\'ern-$3/2$ &  Gaussian Processes \\
    $(1+\sqrt{5}(\frac{r}{\sigma}) + \frac53 (\frac{r}{\sigma})^2)\exp(-\frac{\sqrt{5}r}{\sigma})$ & Mat\'ern-$5/2$&  Gaussian Processes
\end{tabular}
\caption{List of kernels of the form $\kappa(\mat{x},\mat{y}) = \varphi(r)$,  where $r = \|\mat{x}-\mat{y}\|_2$ and $\sigma$ is the scale parameter.} 
\end{table}

\begin{figure}[!ht]
    \centering
    \includegraphics[scale=.35]{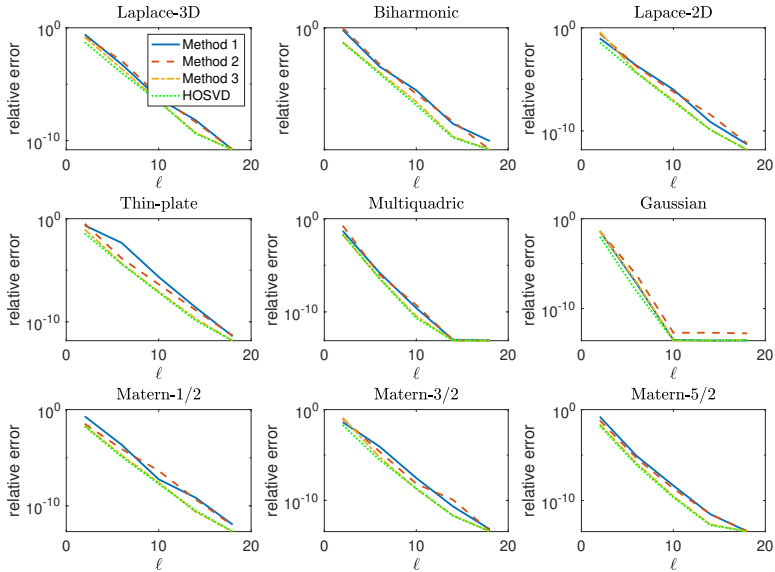}
    \caption{Relative error of Methods 1, 2, and 3 with increasing target rank $r_t$ values compared to HOSVD for nine different kernels.}
    \label{fig:kernelsh}
\end{figure}

We compare relative error with increasing target rank $\ell = r_t + p$. To provide a reasonable comparison in this case, we plot the performance of our algorithms against the relative error of using an HOSVD of $\ten{M}$ instead of the tensor methods we proposed in Section~\ref{sec:lowrankalgs}. These results are plotted in Figure~\ref{fig:kernelsh}. We see that for all five kernels, the accuracy is very similar for each algorithm, with Method 3 and HOSVD as the most accurate.

\begin{figure}[!ht]
    \centering
    \includegraphics[scale=.35]{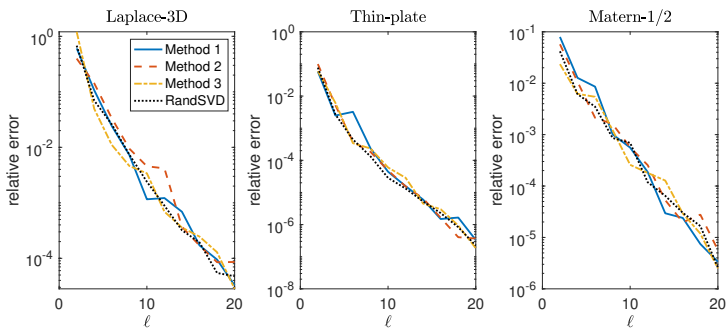}
    \caption{Relative error from Methods 1, 2, and 3 compared to RandSVD with increasing target rank $\ell =r_t + p$. The low-rank approximation for Methods 1, 2, and 3 were recompressed to rank $\ell$ (from rank $\ell^2$).}
    \label{fig:kernels_compress}
\end{figure}

\paragraph{Experiment 2: Comparison with randomized SVD} 
The tensor based methods can be used to produce the low-rank approximation~\eqref{eqn:lowrank2} which has a rank $\ell^D$, however, an alternative approach is to compress the approximate low-rank representation using Chebyshev interpolation $\mat{F}_s\mat{M}\mat{F}_t^\top$ (see~\eqref{eqn:lowrank1}), which has rank $n^D$, using matrix techniques (e.g., SVD or RandSVD). For fair comparison, we compress the output of all the approximations to the same target rank $\ell$ using the approach described in obtaining~\eqref{eqn:lowrank3} (note that the rank of the low-rank approximation before compression is $\ell^2$). The relative error computed for the recompressed tensor-based approximations compared to a RandSVD of the same rank is plotted in Figure~\ref{fig:kernels_compress} with increasing target rank $r_t$. For all five kernels, the accuracy of the tensor based approaches is comparable to that using matrix techniques. However, the tensor-based methods are more advantageous (since Method 2 has fewer kernel evaluations and computational cost, and Method 3 requires generating fewer random numbers).

\begin{figure}[!ht]
    \centering
    \includegraphics[scale=.35]{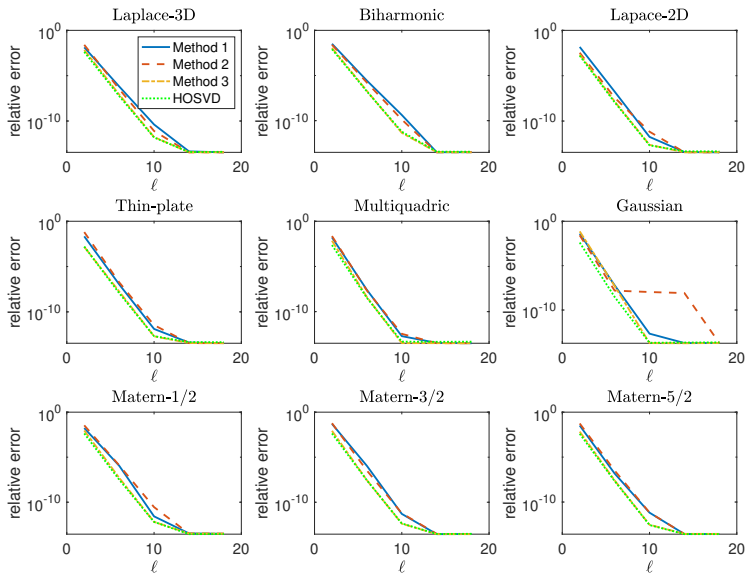}
    \caption{Relative error from Methods 1, 2, 3, and HOSVD with increasing $\ell$ values for the different kernels in three spatial dimensions.}
    \label{fig:3dkernel}
\end{figure}

\paragraph{Experiment 3: Three spatial dimensions}
We now consider experiments in three spatial dimensions. For each of these experiments, we generate the boxes as in the previous subsection. Both the source and target boxes are of size $L=5$ and are placed in the first octant. The source box has one corner at the origin, and the distance between the origin and the corresponding corner of the target box is $D = 15$. In these boxes, we generate $N_s = N_t = 500$ source and target points that are randomly spaced. Furthermore, we take $n = 18$ Chebyshev nodes for each dimension, with block size $n_b = 6$. We plot the relative error from Methods 1, 2, 3, and HOSVD as the target rank $\ell$ increases for five different kernels. These error values are shown in Figure~\ref{fig:3dkernel}. HOSVD and Method 3 are the most accurate for each kernel, but all relative error values are very similar to each other. We see that Method 2 is less accurate for the Gaussian kernel for the larger values of $\ell$.

\subsection{Application to Gaussian Processes}\label{ssec:gauss}
In Gaussian processes involving  large datasets~\cite{liu2020gaussian}, it is desirable to construct a global low-rank approximation when the sources and target points are the same. That is, $\mc{X} = \mc{Y}$, and we want to compute low-rank approximations of the form $\matc{K}(\mc{X},\mc{X})$. Here the implicit assumption is that the kernel is sufficiently smooth so that the matrix admits a low-rank approximation. There are approaches in the literature to generate low-rank approximations such as Nystr\"om~\cite{li2014large}, Structured Kernel Interpolation~\cite{wilson2015kernel}, block-structure~\cite{si2014memory}, etc. 

Since the resulting matrix is symmetric and positive definite, we would like to preserve the symmetry in the low-rank approximation. We can easily adapt the methods proposed in Sections~\ref{sec:lowrankalgs} to handle this. First, by the symmetry of the problem, we have $\mc{B}_s = \mc{B}_t$ and, therefore, $c_j = a_j$ and $d_j = b_j$ for $j=1,\dots,D$. Furthermore, we compute the factor matrices only along the first $D$ modes. This gives us the compressed tensor representation $\ten{M} \approx [\ten{G}; \mat{A}_1,\dots,\mat{A}_D, \mat{A}_1,\dots\,\mat{A}_D]$ with target rank $(r_t,\dots,r_t)$ which we can also express as
\[ \ten{M} \approx \ten{G} \left(\bigtimes_{j=1}^D \mat{A}_j \right) \left(\bigtimes_{j=D+1}^{2D} \mat{A}_{j-D}\right), \qquad \ten{G} = \ten{M} \left(\bigtimes_{j=1}^D \mat{A}_j^\top \right) \left(\bigtimes_{j=D+1}^{2D} \mat{A}_{j-D}^\top\right) . \]
Using this tensor approximation, we can obtain a low-rank approximation kernel approximation
\[\matc{K}(\mc{X},\mc{X}) \approx \mathat{F}\mathat{M} \mathat{F}^\top, \]
where $\mathat{F} = \mathat{U}_D \ltimes \mathat{U}_{D-1} \ltimes \dots \ltimes \mathat{U}_1$. The storage cost of the low-rank approximation is $2r_tDN_s + r_t^{2D}$; if $\mathat{F}$ is formed is explicitly, then the cost of storing is $2r_t^DN_s + r_t^{2D}$ entries. For simplicity, we have assumed the amount of oversampling is $p =0$, but it is straightforward to incorporate additional oversampling. 

We use the daily precipitation data set corresponding to $5,500$ weather stations in the US for the year 2010~\cite{datawebsite}. The details of this test problem are given in~\cite{dong2017scalable}. In total, this data set has $N_s = 628,474$ points. Our goal is only to compare the accuracy of the kernel matrix approximations and not to study its application to Gaussian Process regression. To this end, we denote $\mat{K} = \matc{K}(\mc{X},\mc{X})$ and $\mathat{K} = \mathat{F}\mathat{M} \mathat{F}^\top $. We compute the relative error as 
\[ \text{rel err} = \frac{|\text{trace}(\mat{K})- \text{trace}(\mathat{K})|}{\text{trace}(\mat{K})}. \]
We use this error measure since it only requires computing the diagonals of $\mat{K}$ and $\mathat{K}$. 
We take the kernel to be 
\[ \kappa(\mat{x},\mat{y}) = \exp \left( -\sum_{j=1}^3\frac{(y_j-x_j)^2 }{\sigma_j^2}\right)^{1/2}.\]
The scale parameters are taken to be $\sigma_1 = 0.8\cdot 80$, $\sigma_2 =  0.8\cdot30$, $\sigma_3 =  0.8\cdot 365.$
We take the number of Chebyshev points $n=27$ and the number of blocks $n_b=9$ and the results are reported in Table~\ref{tab:gpdataset}. Once again, we see that all three methods are highly accurate and the accuracy of the methods are comparable to that using HOSVD. In the worst case, Method 2 only required $\approx 1.3\%$ of the total function evaluations.

\begin{table}[!ht]
    \centering
    \begin{tabular}{c|c|c|c|c}
         & $r_t = 2$ & $r_t = 4$ & $r_t = 6$ &  $r_t = 8$ \\ \hline
         HOSVD & $8.06\times 10^{-2}$ & $2.49\times 10^{-4}$ & $3.24\times 10^{-7}$ & $2.14\times 10^{-10}$ \\ 
         Method 1    & $4.04\times 10^{-2}$ & $5.39\times 10^{-4}$ & $1.40\times 10^{-6}$ & $9.84\times 10^{-10}$ \\
    Method 2     & $1.31\times 10^{-1}$ & $5.13\times 10^{-4}$ & $3.57\times 10^{-6}$ & $1.28\times 10^{-9}$ 
    \end{tabular}
    \caption{Accuracy of the low-rank approximation corresponding to the precipitation data set. }
    \label{tab:gpdataset}
\end{table}

\section{Conclusions and future work}
In summary, we proposed a tensor-based function approximation, developed new randomized algorithms for tensor compression in the Tucker format, and demonstrated an application of these approaches to constructing low-rank approximations to kernel matrices. When the function is sufficiently smooth, the resulting tensor-based approximations are accurate and storage efficient and can be used as a surrogate model in a variety of applications. The kernel low-rank approximations are efficient to compute and store; both the computational and storage costs are competitive with other methods.  An important future direction that we have not addressed in this paper is to provide an adaptive approach for computing the function approximation to a desired accuracy. This could, for example, use the adaptive techniques developed in~\cite{hashemi2017chebfun,dolgov2021functional} that uses multiple computational phases. Another direction for future work is to consider tensor compression techniques using the Tensor Train format rather than the Tucker format and compare it with related techniques~\cite{gorodetsky2019continuous,bigoni2016spectral}. This has the potential to be effective for very high-dimensional problems.

 \appendix
 \section{Proof of Theorem~\ref{thm:approx}}\label{sec:appenda}
For each $j = 1,\dots,N$, let $\mat{\Pi}_j = \mat{A}_j\mat{P}_j^\top$.
First recall that $\tenh{M} = \ten{M} \bigtimes_{j=1}^N \mat{\Pi}_j$. Then considering the term $\ten{M}-\tenh{M}$, we can add and subtract the terms $\ten{M} \bigtimes_{i=1}^{j} \mat{\Pi}_i$ for $j = 1,\dots,N-1$ as follows.
\begin{equation*}
    \begin{aligned}
        \ten{M}-\tenh{M} &= \ten{M}-\ten{M} \bigtimes_{j=1}^N \mat{\Pi}_j \\
        &= \ten{M} \times_1 (\mat{I}-\mat{\Pi}_1) + \ten{M} \times_1 \mat{\Pi}_1 \times_2 (\mat{I}-\mat{\Pi}_2) + \dots + \ten{M} \bigtimes_{j=1}^{N-1} \mat{\Pi}_j \times_N (\mat{I}-\mat{\Pi}_N).
    \end{aligned}
\end{equation*}
Then, taking the Frobenius norm and applying the triangle inequality, we have 
\begin{equation*}
    \|\ten{M}-\tenh{M}\|_F \leq \sum_{j=1}^N \| \ten{M} \bigtimes_{i=1}^{j-1} \mat{\Pi}_i \times_j (\mat{I}-\mat{\Pi}_j) \|_F.
\end{equation*}
Using the linearity of expectations we get 
\begin{equation}\label{eq:sum_exp}
    \mb{E}_{\{\mat{\Omega}\}_{k=1}^N} \|\ten{M}-\tenh{M}\|_F \leq \sum_{j=1}^N \mb{E}_{\mat{\Omega}_j} \| \ten{M} \bigtimes_{i=1}^{j-1} \mat{\Pi}_i \times_j (\mat{I}-\mat{\Pi}_j) \|_F.
\end{equation}
Now consider the term $\|\ten{M} \bigtimes_{i=1}^{j-1} \mat{\Pi}_i \times_j (\mat{I}-\mat{\Pi}_j) \|_F$. We can use the submultiplicativity property $\| \mat{A}\mat{B}\|_F \leq \| \mat{A}\|_2 \|\mat{B}\|_F$ to obtain
\begin{equation}\label{eq:pi}
    \|\ten{M} \bigtimes_{i=1}^{j-1} \mat{\Pi}_i \times_j (\mat{I}-\mat{\Pi}_j) \|_F \leq \prod_{i=1}^{j-1} \|\mat{\Pi}_i\|_2 \| \ten{M} \times_j (\mat{I}-\mat{\Pi}_j)\|_F.
\end{equation}
Note that $\mat{\Pi}_j\mat{Q}_j \mat{Q}_j^\top = \mat{Q}_j \mat{Q}_j^\top$, meaning that $\mat{I}-\mat{\Pi}_j = (\mat{I}-\mat{\Pi}_j)(\mat{I}-\mat{Q}_j\mat{Q}_j^\top)$. We can then rewrite each term $\|\ten{M} \times_j (\mat{I}-\mat{\Pi}_j)\|_F$ as 
\begin{equation}\label{eq:submult}
\begin{aligned}
    \| \ten{M} \times_j (\mat{I}-\mat{\Pi}_j)\|_F &= \| \ten{M} \times_j (\mat{I}-\mat{\Pi}_j)(\mat{I}-\mat{Q}_j \mat{Q}_j^\top)\|_F \\
    &\leq \|\ten{M} \times_j (\mat{I}-\mat{Q}_j\mat{Q}_j^\top)\|_F \|\mat{I}-\mat{\Pi}_j\|_2.
\end{aligned}
\end{equation}
%Recall
We note that $\mat{\Pi}_j \neq \mat{I}$ since $\rank(\mat{\Pi}_j) \leq \rank(\mat{Q}_j) = r_t + p < n$, and $\mat{\Pi}_j \neq \mathbf{0}$ since $\mat{Q}_j$ has orthonormal columns  and $\mat{P}_j^\top \mat{Q}_j$ is invertible. Therefore, using the main result of~\cite{szyld2006many}, $\|\mat{I} - \mat{\Pi}_j\|_2 = \|\mat{\Pi}_j\|_2$. Then using \cite[Lemma 2.1]{drmac2018discrete}, 
$\|\mat{I}-\mat{\Pi}_j\|_2 = \|\mat{\Pi}_j\|_2 = \|(\mat{P}_j^\top\mat{Q}_j)^{-1}\|_2 \leq g(n,\ell).$
Combining this result with equations \eqref{eq:pi} and \eqref{eq:submult}, we obtain
\begin{equation*}
    \begin{aligned}
    \|\ten{M} \bigtimes_{i=1}^{j-1} \mat{\Pi}_j \times (\mat{I}-\mat{\Pi}_j)\|_F &\leq  \left(g(n,\ell)\right)^j \|\ten{M} \times_j (\mat{I}-\mat{Q}_j\mat{Q}_j^\top)\|_F \\
    &= \left(g(n,\ell)\right)^j \|(\mat{I}-\mat{Q}_j\mat{Q}_j^\top)\mat{M}_{(j)}\|_F.
    \end{aligned}
\end{equation*}
The last equality comes by unfolding the term $\ten{M} \times_j (\mat{I}-\mat{Q}_j\mat{Q}_j^\top)$ along mode-$j$. 
By taking expectations and applying \cite[Theorem 3]{zhang2018randomized}, we have
\begin{equation*}
\begin{aligned}
    \mb{E}_{\mat{\Omega}_j}\|\ten{M} \bigtimes_{i=1}^{j-1} \mat{\Pi}_j \times_j (\mat{I}-\mat{\Pi}_j)\|_F &\leq (g(n,\ell))^j \mb{E}_{\mat{\Omega}_j} \|(\mat{I}-\mat{Q}_j\mat{Q}_j^\top)\mat{M}_{(j)}\|_F \\
    &\leq (g(n,\ell))^j \left(\left(1+\frac{r_t}{p-1} \right) \sum_{i=r_t+1}^n \sigma_i^2\left(\mat{M}_{(j)}\right)\right)^{1/2}.
    \end{aligned}
\end{equation*}
Combining this result with \eqref{eq:sum_exp}, we obtain the desired bound.

\bibliographystyle{abbrv}
\bibliography{ref}

\begin{thebibliography}{10}

\bibitem{datawebsite}
{U.S.} hourly precipitation data,
  https://catalog.data.gov/dataset/u-s-hourly-precipitation-data.

\bibitem{barnes1986hierarchical}
J.~Barnes and P.~Hut.
\newblock A hierarchical $\mc{O}(n \log n)$ force-calculation algorithm.
\newblock {\em Nature}, 324(6096):446--449, 1986.

\bibitem{bebendorf2009recompression}
M.~Bebendorf and S.~Kunis.
\newblock Recompression techniques for adaptive cross approximation.
\newblock {\em The Journal of Integral Equations and Applications}, pages
  331--357, 2009.

\bibitem{bigoni2016spectral}
D.~Bigoni, A.~P. Engsig-Karup, and Y.~M. Marzouk.
\newblock Spectral tensor-train decomposition.
\newblock {\em SIAM Journal on Scientific Computing}, 38(4):A2405--A2439, 2016.

\bibitem{sfuwebsite}
D.~Bingham.
\newblock Virtual library of simulation experiments: {T}est functions and
  datasets, {http://www.sfu.ca/~ssurjano/index.html}.

\bibitem{borm2003hierarchical}
S.~B{\"o}rm, L.~Grasedyck, and W.~Hackbusch.
\newblock Hierarchical matrices.
\newblock {\em Lecture notes}, 21:2003, 2003.

\bibitem{bungartz_griebel_2004}
H.-J. Bungartz and M.~Griebel.
\newblock Sparse grids.
\newblock {\em Acta Numerica}, 13:147–269, 2004.

\bibitem{cambier2019fast}
L.~Cambier and E.~Darve.
\newblock Fast low-rank kernel matrix factorization using skeletonized
  interpolation.
\newblock {\em SIAM Journal on Scientific Computing}, 41(3):A1652--A1680, 2019.

\bibitem{chen2018fast}
C.~Chen, S.~Aubry, T.~Oppelstrup, A.~Arsenlis, and E.~Darve.
\newblock Fast algorithms for evaluating the stress field of dislocation lines
  in anisotropic elastic media.
\newblock {\em Modelling and Simulation in Materials Science and Engineering},
  26(4):045007, 2018.

\bibitem{constantine2017global}
P.~G. Constantine and P.~Diaz.
\newblock Global sensitivity metrics from active subspaces.
\newblock {\em Reliability Engineering \& System Safety}, 162:1--13, 2017.

\bibitem{de2000multilinear}
L.~De~Lathauwer, B.~De~Moor, and J.~Vandewalle.
\newblock A multilinear singular value decomposition.
\newblock {\em SIAM journal on Matrix Analysis and Applications},
  21(4):1253--1278, 2000.

\bibitem{dolgov2021functional}
S.~Dolgov, D.~Kressner, and C.~Str{\"o}ssner.
\newblock Functional {T}ucker approximation using {C}hebyshev interpolation.
\newblock {\em SIAM Journal on Scientific Computing}, 43(3):A2190--A2210, 2021.

\bibitem{dong2017scalable}
K.~Dong, D.~Eriksson, H.~Nickisch, D.~Bindel, and A.~G. Wilson.
\newblock Scalable log determinants for gaussian process kernel learning.
\newblock In I.~Guyon, U.~V. Luxburg, S.~Bengio, H.~Wallach, R.~Fergus,
  S.~Vishwanathan, and R.~Garnett, editors, {\em Advances in Neural Information
  Processing Systems}, volume~30. Curran Associates, Inc., 2017.

\bibitem{drmac2018discrete}
Z.~Drma\v{c} and A.~K. Saibaba.
\newblock The discrete empirical interpolation method: Canonical structure and
  formulation in weighted inner product spaces.
\newblock {\em SIAM Journal on Matrix Analysis and Applications},
  39(3):1152--1180, 2018.

\bibitem{fong2009black}
W.~Fong and E.~Darve.
\newblock The black-box fast multipole method.
\newblock {\em Journal of Computational Physics}, 228(23):8712--8725, 2009.

\bibitem{gass2018chebyshev}
M.~Ga{\ss}, K.~Glau, M.~Mahlstedt, and M.~Mair.
\newblock {C}hebyshev interpolation for parametric option pricing.
\newblock {\em Finance and Stochastics}, 22(3):701--731, 2018.

\bibitem{golub2013matrix}
G.~H. Golub and C.~F. Van~Loan.
\newblock {\em Matrix computations}.
\newblock Johns Hopkins Studies in the Mathematical Sciences. Johns Hopkins
  University Press, Baltimore, MD, fourth edition, 2013.

\bibitem{gorodetsky2019continuous}
A.~Gorodetsky, S.~Karaman, and Y.~Marzouk.
\newblock A continuous analogue of the tensor-train decomposition.
\newblock {\em Computer Methods in Applied Mechanics and Engineering},
  347:59--84, 2019.

\bibitem{grasedyck2003construction}
L.~Grasedyck and W.~Hackbusch.
\newblock Construction and arithmetics of $\mc{H}$-matrices.
\newblock {\em Computing}, 70(4):295--334, 2003.

\bibitem{greengard1987fast}
L.~Greengard and V.~Rokhlin.
\newblock A fast algorithm for particle simulations.
\newblock {\em Journal of computational physics}, 73(2):325--348, 1987.

\bibitem{greengard1991fast}
L.~Greengard and J.~Strain.
\newblock The fast {G}auss transform.
\newblock {\em SIAM Journal on Scientific and Statistical Computing},
  12(1):79--94, 1991.

\bibitem{griebel2019analysis}
M.~Griebel and H.~Harbrecht.
\newblock Analysis of tensor approximation schemes for continuous functions.
\newblock {\em arXiv preprint arXiv:1903.04234}, 2019.

\bibitem{gu1996efficient}
M.~Gu and S.~C. Eisenstat.
\newblock Efficient algorithms for computing a strong rank-revealing {QR}
  factorization.
\newblock {\em SIAM Journal on Scientific Computing}, 17(4):848--869, 1996.

\bibitem{hackbusch1999sparse}
W.~Hackbusch.
\newblock A sparse matrix arithmetic based on $\mc{H}$-matrices. part i:
  Introduction to $\mc{H}$-matrices.
\newblock {\em Computing}, 62(2):89--108, 1999.

\bibitem{hackbusch2015hierarchical}
W.~Hackbusch.
\newblock {\em Hierarchical matrices: algorithms and analysis}, volume~49.
\newblock Springer, 2015.

\bibitem{hackbusch2002data}
W.~Hackbusch and S.~B{\"o}rm.
\newblock Data-sparse approximation by adaptive $\mc{H}^2$-matrices.
\newblock {\em Computing}, 69(1):1--35, 2002.

\bibitem{halko2011finding}
N.~Halko, P.-G. Martinsson, and J.~A. Tropp.
\newblock Finding structure with randomness: Probabilistic algorithms for
  constructing approximate matrix decompositions.
\newblock {\em SIAM review}, 53(2):217--288, 2011.

\bibitem{hashemi2017chebfun}
B.~Hashemi and L.~N. Trefethen.
\newblock Chebfun in three dimensions.
\newblock {\em SIAM Journal on Scientific Computing}, 39(5):C341--C363, 2017.

\bibitem{kolda2009tensor}
T.~G. Kolda and B.~W. Bader.
\newblock Tensor decompositions and applications.
\newblock {\em SIAM review}, 51(3):455--500, 2009.

\bibitem{li2014large}
M.~Li, W.~Bi, J.~T. Kwok, and B.-L. Lu.
\newblock Large-scale {N}ystr{\"o}m kernel matrix approximation using
  randomized {SVD}.
\newblock {\em IEEE transactions on neural networks and learning systems},
  26(1):152--164, 2014.

\bibitem{liu2020gaussian}
H.~Liu, Y.-S. Ong, X.~Shen, and J.~Cai.
\newblock When {G}aussian process meets big data: A review of scalable {GP}s.
\newblock {\em IEEE transactions on neural networks and learning systems},
  31(11):4405--4423, 2020.

\bibitem{mason2002chebyshev}
J.~C. Mason and D.~C. Handscomb.
\newblock {\em Chebyshev polynomials}.
\newblock CRC press, 2002.

\bibitem{minster2020randomized}
R.~Minster, A.~K. Saibaba, and M.~E. Kilmer.
\newblock Randomized algorithms for low-rank tensor decompositions in the
  {T}ucker format.
\newblock {\em SIAM Journal on Mathematics of Data Science}, 2(1):189--215,
  2020.

\bibitem{rai2019randomized}
P.~Rai, H.~Kolla, L.~Cannada, and A.~Gorodetsky.
\newblock Randomized functional sparse {T}ucker tensor for compression and fast
  visualization of scientific data.
\newblock {\em arXiv preprint arXiv:1907.05884}, 2019.

\bibitem{sauter2010boundary}
S.~A. Sauter and C.~Schwab.
\newblock Boundary element methods.
\newblock In {\em Boundary Element Methods}, pages 183--287. Springer, 2010.

\bibitem{si2014memory}
S.~Si, C.-J. Hsieh, and I.~Dhillon.
\newblock Memory efficient kernel approximation.
\newblock In {\em International Conference on Machine Learning}, pages
  701--709. PMLR, 2014.

\bibitem{szyld2006many}
D.~B. Szyld.
\newblock The many proofs of an identity on the norm of oblique projections.
\newblock {\em Numerical Algorithms}, 42(3-4):309--323, 2006.

\bibitem{trefethen2017multivariate}
L.~Trefethen.
\newblock Multivariate polynomial approximation in the hypercube.
\newblock {\em Proceedings of the American Mathematical Society},
  145(11):4837--4844, 2017.

\bibitem{trefethen2013approximation}
L.~N. Trefethen.
\newblock {\em Approximation theory and approximation practice}.
\newblock Society for Industrial and Applied Mathematics (SIAM), Philadelphia,
  PA, 2013.

\bibitem{wilson2015kernel}
A.~Wilson and H.~Nickisch.
\newblock Kernel interpolation for scalable structured {G}aussian processes
  ({KISS-GP}).
\newblock In {\em International Conference on Machine Learning}, pages
  1775--1784. PMLR, 2015.

\bibitem{xu2016chebyshev}
K.~Xu.
\newblock The {C}hebyshev points of the first kind.
\newblock {\em Applied Numerical Mathematics}, 102:17--30, 2016.

\bibitem{xu2018low}
Z.~Xu, L.~Cambier, F.-H. Rouet, P.~L'Eplatennier, Y.~Huang, C.~Ashcraft, and
  E.~Darve.
\newblock Low-rank kernel matrix approximation using skeletonized interpolation
  with endo-or exo-vertices.
\newblock {\em arXiv preprint arXiv:1807.04787}, 2018.

\bibitem{ye2020analytical}
X.~Ye, J.~Xia, and L.~Ying.
\newblock Analytical low-rank compression via proxy point selection.
\newblock {\em SIAM Journal on Matrix Analysis and Applications},
  41(3):1059--1085, 2020.

\bibitem{ying2004kernel}
L.~Ying, G.~Biros, and D.~Zorin.
\newblock A kernel-independent adaptive fast multipole algorithm in two and
  three dimensions.
\newblock {\em Journal of Computational Physics}, 196(2):591--626, 2004.

\bibitem{zhang2018randomized}
J.~Zhang, A.~K. Saibaba, M.~E. Kilmer, and S.~Aeron.
\newblock A randomized tensor singular value decomposition based on the
  t-product.
\newblock {\em Numerical Linear Algebra with Applications}, 25(5):e2179, 2018.

\end{thebibliography}
\end{document}